\documentclass[smallextended]{svjour3}

\usepackage{amsmath}
\usepackage{amssymb}
\usepackage{bm}
\usepackage{enumitem}
\usepackage{etoolbox}
\usepackage{filecontents}
\usepackage{hyperref}
\usepackage{natbib}
\usepackage{nomencl}
\usepackage{mathtools}
\usepackage{pgfplots}
\usepackage{pgfplotstable}
\usepackage{subcaption}

\hypersetup{colorlinks=true,linkcolor=black,citecolor=black,filecolor=black,urlcolor=black}
\usetikzlibrary{pgfplots.statistics,pgfplots.colorbrewer,patterns}
\pgfplotsset{compat=1.14,cycle list/Set1-9,
        /pgfplots/ybar legend/.style={
        /pgfplots/legend image code/.code={%
        \draw[##1,/tikz/.cd,bar width=3pt,yshift=-0.2em,bar shift=0pt]
                plot coordinates {(0cm,0.8em)};},
},}

\pgfdeclarepatternformonly{finedots}{\pgfqpoint{-0.5pt}{-0.5pt}}{\pgfqpoint{0pt}{0pt}}{\pgfqpoint{0.5pt}{0.5pt}}%
{
    \pgfpathcircle{\pgfqpoint{0pt}{0pt}}{.2pt}
    \pgfpathcircle{\pgfqpoint{1pt}{1pt}}{.2pt}
    \pgfusepath{fill}
}

\begin{document}

\title{Natural Gas Maximal Load Delivery for Multi-contingency Analysis}
\author{Byron Tasseff \and Carleton Coffrin \and Russell Bent \and  Kaarthik Sundar \and Anatoly Zlotnik}
    \institute{
    \href{https://orcid.org/0000-0002-5043-8305}{\includegraphics[width=5.5px]{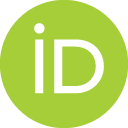}} B. Tasseff \at
    Los Alamos National Laboratory, Los Alamos, New Mexico, USA \\
    University of Michigan, Ann Arbor, Michigan, USA \\
    \email{btasseff@lanl.gov} \and
    \href{https://orcid.org/0000-0003-3238-1699}{\includegraphics[width=5.5px]{image/orcid.png}} C. Coffrin,
    \href{https://orcid.org/0000-0002-7300-151X}{\includegraphics[width=5.5px]{image/orcid.png}} R. Bent,
    \href{https://orcid.org/0000-0002-6928-449X}{\includegraphics[width=5.5px]{image/orcid.png}} K. Sundar, and
    \href{https://orcid.org/0000-0002-2646-8264}{\includegraphics[width=5.5px]{image/orcid.png}} A. Zlotnik \at
    Los Alamos National Laboratory, Los Alamos, New Mexico, USA}
\maketitle

\begin{abstract}
As the use of renewable generation has increased, electric power systems have become increasingly reliant on natural gas-fired power plants as fast ramping sources for meeting fluctuating bulk power demands.
This dependence has introduced new vulnerabilities to the power grid, including disruptions to gas transmission networks from natural and man-made disasters.
To address the operational challenges arising from these disruptions, we consider the task of determining a feasible steady-state operating point for a damaged gas pipeline network while ensuring the maximal delivery of load.
We formulate the mixed-integer nonconvex maximal load delivery (MLD) problem, which proves difficult to solve on large-scale networks.
To address this challenge, we present a mixed-integer convex relaxation of the MLD problem and use it to determine bounds on the transport capacity of a gas pipeline system.
To demonstrate the effectiveness of the relaxation, the exact and relaxed formulations are compared across a large number of randomized damage scenarios on nine natural gas pipeline network models ranging in size from $11$ to $4197$ junctions.
A proof of concept application, which assumes network damage from a set of synthetically generated earthquakes, is also presented to demonstrate the utility of the proposed optimization-based capacity evaluation in the context of risk assessment for natural disasters.
For all but the largest network, the relaxation-based method is found to be suitable for use in evaluating the impacts of multi-contingency network disruptions, often converging to the optimal solution of the relaxed formulation in less than ten seconds.
\end{abstract}


\section{Introduction}
\label{section:introduction}
Globally, electricity generation capacity is expected to grow from $21.6$ trillion kilowatt-hours (kWh) in 2012 to $36.5$ trillion kWh in 2040.
Of this capacity, generation from natural gas is expected to grow from $22\%$ to $28\%$ \citep{conti2016international}.
This growing dependence implies the susceptibility of power systems to upstream disruptions in gas transmission networks.
Historical examples include the 2014 polar vortex, where natural gas delivery curtailments led to approximately $25\%$ of the total generation outages in the PJM interconnection \citep{pjm2014}.
A separate example is the 2014 South Napa earthquake, which highlighted the vulnerability of gas networks to ground failure \citep{johnson2016mw}.
Other hazards, including floods and wildfires, can also have severe impacts on gas networks \citep{infrastructure2014united}.
Aside from such downstream effects on the power grid, these disruptions further inhibit the transport of fuel for residential heating, which provides essential temperature control to many individual homes during winter months.
Mitigating the effects of these disruptions is critical to the resilience of gas and power delivery networks.
To that end, this study examines how to compute the optimal response to a large-scale multi-contingency gas pipeline network disruption, whose origin (e.g., a natural hazard or sophisticated attack) is treated agnostically.

\begin{figure}[t]
    \begin{center}
        \includegraphics[width=\textwidth]{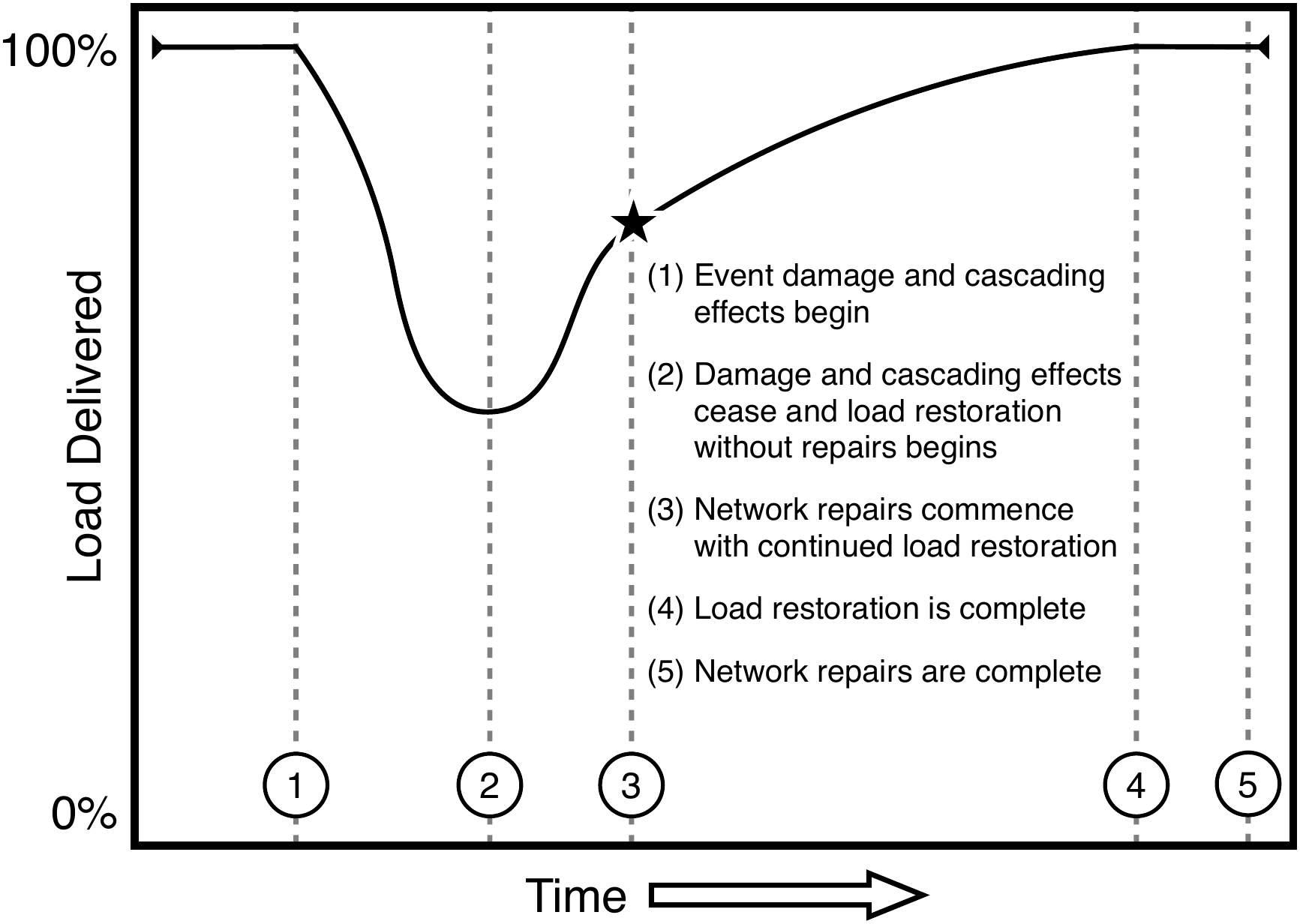}
    \end{center}
    \caption{Illustration of a gas network's response to an extreme event, similar to the power network response in \citet{coffrin2018relaxations}. The star indicates the point in the damage and restoration process that is examined in this study with an optimization-based assessment.}
    \label{figure:restoration}
\end{figure}

The scope of the response measures considered in this study is illustrated at a high level in Figure \ref{figure:restoration}.
When a hazard event begins, (1) delivery of load decreases as gas network components are damaged and cascading effects begin, and (2) cascading effects cease, and a new stable operating point is found.
After (2), some amount of load can be restored through operational methods until (3) network repairs commence.
These repairs are conducted in accord with other restoration processes until (4) all load can be delivered.
Repairs continue until (5) all network components are operational.
Addressing the complete scope of Figure \ref{figure:restoration} is a substantial and complex task.
In this study, we focus on determining optimal steady-state operating points between events of types (2) and (4), i.e., operational restoration decisions that enable the delivery of a maximal amount of load in a damaged gas pipeline network.

Several commercial tools exist for analyzing the operation of gas pipelines in the steady-state and transient regimes, including the \textsc{NextGen} pipeline simulation suite \citep{nextgen}, Energy Solutions' gas management software \citep{energysolutions}, and \textsc{Atmos Pipe} \citep{atmos}.
These tools are designed for capacity planners to simulate the operation of a gas pipeline network under various physical conditions and flow nominations rather than for resilience analysis.
When one or more components are set to be nonoperational (e.g., to simulate damage caused by a natural disaster) and must be removed from the network model, this must typically be done manually in such commercially available software tools.
Therefore, evaluating the impact of many thousands of multi-component outage scenarios using existing tools is too labor-intensive to be practical.
Furthermore, because some deliveries (demand) and receipts (supply) may need to be adjusted (or omitted) because of component outages, a new feasible operating solution, including flow allocation and compressor settings, needs to be quickly obtained for each scenario.  It is therefore unclear how existing commercial tools can be used for probabilistic risk assessment in this setting.
This motivates a mathematical approach that can determine a feasible natural gas pipeline operating point that maximizes delivery subject to the multi-contingency considerations.

In this study, we formalize this task as the steady-state Maximal Load Delivery (MLD) problem.
Informally, the problem can be stated as follows: given a severely damaged gas pipeline network in which a number of components have become nonoperational, we seek to maximize the amount of prioritized load that can be served in the damaged network subject to steady-state pipeline physical flow, capacity limits, pressure bounds, and other operating requirements.
The nonlinear physics of gas transport and the discrete nature of operations in the network (i.e., the opening and closing of valves) make this problem a challenging mixed-integer nonconvex program.
To address this, we introduce a mixed-integer convex relaxation of the problem, which is found to be a reliable means for determining bounds on the maximum deliverable load.

Recent interest in large-scale gas network planning and control has led to a wide variety of optimization-based applications and methodologies.
\citet{HILLER2018797} provide a summary of studies connected to the optimization-based evaluation of gas network capacities.
\citet{schmidt2016high} present detailed steady-state models and approximations of network components for use in optimization applications.
\citet{hahn2017mixed}, \citet{gugat2018mip}, and \citet{hoppmann2019optimal} present mixed-integer programming approaches for optimization problems involving gas transport in the transient regime.
\citet{doi:10.1002/oca.2315,Hante2020} describe relaxation methods for similar mixed-integer control problems.

The methodology of this study is primarily inspired by the success of approaches developed for multi-contingency analysis of power transmission networks.
A similar optimization study on damaged power grids was motivated by the analysis of natural disaster vulnerabilities \citep{coffrin2018relaxations}.
In that study, the MLD problem was formulated for alternating current (AC) power networks, and convex relaxations were developed that allow for the problem's efficient solution on large-scale instances.
Another study introduced convex relaxations for gas network expansion planning \citep{borraz2016convex}.
Similar convex relaxations for gas pipeline flow have been explored subsequently \citep{wu2017adaptive,chen2018steady}.
These relaxations largely inspire the convex reformulation of the MLD problem described in this study.
Extending these studies, we present nonconvex models and convex relaxations for components that were previously not explicitly considered, e.g., resistors.

A number of studies consider problems as for the one developed here.
\citet{sundar2018probabilistic} and \citet{ahumada2021nk} examine the problems of identifying the $k$ components of a power transmission network and gas pipeline network, respectively, whose simultaneous failure maximizes disruption to the network.
Their studies use MLD problems as inner portions of broader bilevel interdiction problems.
In contrast, we seek to determine an optimal operating point for a given disruption rather than find the worst-case disruption.
\citet{bent2018joint} consider the joint expansion planning problem for gas and power networks.
Their study, which links power generation to gas delivery, serves as a reference for future extensions of the methods presented in this study.

The contributions of this study include
\begin{itemize}[noitemsep,topsep=0pt]
    \item The first formulation of the steady-state MLD problem for gas networks;
    \item A mixed-integer nonconvex quadratic reformulation of the MLD problem;
    \item A tractable mixed-integer convex relaxation of the MLD problem;
    \item A rigorous benchmarking of the formulations on instances of various sizes;
    \item A proof of concept MLD analysis for spatially distributed natural hazards.
\end{itemize}

The remaining sections proceed as follows: Section \ref{section:background} formulates the requirements for gas pipeline operational feasibility as a mixed-integer nonconvex program;
Section \ref{section:mld-formulations} formulates the MLD problem as a mixed-integer nonconvex program, then proposes a mixed-integer convex relaxation of the problem;
Section \ref{section:computational_experiments} rigorously benchmarks the formulations across several gas network data sets of various sizes and examines the use of the MLD method for probabilistic risk assessment for natural disasters;
and Section \ref{section:concluding_remarks} concludes the paper.

\section{Background}
\label{section:background}
This section presents the background for modeling natural gas transmission networks in a steady-state regime.
In Section \ref{subsection:pipeline_modeling_assumptions}, the transient dynamics for gas pipelines are stated, and standard physical assumptions are imposed to derive the Weymouth equation for steady-state turbulent flow.
This nonconvex equation describes the loss in potential (i.e., pressure) along each pipe as a nonlinear function of flow and represents one of the primary modeling challenges addressed in subsequent sections.
Then, in Section \ref{subsection:network-modeling}, additional physical constraints required for feasible gas network operation are described.

\subsection{Pipeline Modeling Assumptions}
\label{subsection:pipeline_modeling_assumptions}
The dynamics of gas pipeline flow are accurately modeled via the following set of partial differential equations (PDEs) over space and time \citep{osiadacz1987simulation}:
\begin{subequations}
    \begin{align}
        \frac{\partial \phi}{\partial x} + \frac{\partial \rho}{\partial t} &= 0, \label{eqn:pde-1} \\
        \frac{\partial (\phi v)}{\partial x} + \frac{\partial \phi}{\partial t} + \frac{\partial p}{\partial x} &= -\frac{\lambda}{2 D} \phi \lvert v \rvert - g \sin \alpha \rho, \label{eqn:pde-2} \\
        p - z R T \rho &= 0 \label{eqn:pde-3}.
    \end{align}\label{eqn:pde}%
\end{subequations}
Here, Equation \eqref{eqn:pde-1} represents conservation of mass, Equation \eqref{eqn:pde-2} represents conservation of momentum, and Equation \eqref{eqn:pde-3} is the equation of state that relates the pressure and density.
In these PDEs, the gas mass flux $\phi$, density $\rho$, velocity $v$, and pressure $p$ evolve over both space and time, $x$ and $t$, respectively.
In Equation \eqref{eqn:pde-2}, $\lambda$ denotes the friction factor of the pipe and, unless given, is computed via an approximation of the Colebrook-White equation for turbulent flow.
That approximation is empirically derived as \citep{zeghadnia2019explicit}
\begin{equation}
    \lambda = \left[2 \log \left(\frac{3.7 D}{\epsilon}\right)\right]^{-2},
\end{equation}
where $D$ denotes the diameter of the pipe wall and $\epsilon$ denotes the absolute pipe roughness (in units of length).
Additionally, in Equations \eqref{eqn:pde-2} and \eqref{eqn:pde-3}, $g$ denotes acceleration due to gravity, $\alpha$ denotes the angle of inclination of the pipe with respect to the horizontal, $z$ denotes the gas compressibility factor, $T$ denotes the temperature of the gas, and $R$ denotes the universal gas constant.
While in practice, the gas compressibility $z$ depends significantly on pressure and temperature in the regime of high pressure gas pipeline flow, in this study, we use an ideal gas equation of state and suppose $z$ to be fixed, as commonly done to explore new concepts in academic studies \citep{schmidt2016high}.

In our study, we assume all pipes are level, which allows the gravitational term of Equation \eqref{eqn:pde-2} to be ignored.
It is also assumed that the temperature $T$ is constant along a pipe.
Finally, we assume that the system has reached a steady state, and thus all time derivatives vanish.
With these assumptions, Equation \eqref{eqn:pde-1} reduces to $\partial_{x} \phi = 0$.
Equation \eqref{eqn:pde-2}, after imposing these assumptions and multiplying by $\rho$, then reduces to
\begin{equation}
    p \frac{\partial p}{\partial x} = -\frac{z R \lambda T}{2 D} \phi \lvert \phi \rvert \implies \frac{\partial p^{2}}{\partial x} = -\frac{z R \lambda T}{D} \phi \lvert \phi \rvert \label{eqn:pde-approx-2}.
\end{equation}
Finally, integrating Equation \eqref{eqn:pde-approx-2} from zero to the length of the pipe, $L$, gives
\begin{equation}
    [p(0)]^{2} - [p(L)]^{2} = \gamma \phi \lvert \phi \rvert \implies [p(L)]^{2} - [p(0)]^{2} = -\gamma \phi \lvert \phi \rvert \label{eqn:weymouth},
\end{equation}
where the resistance term is defined via the relation $\gamma := \frac{z L R \lambda T}{D}$.
Equation \eqref{eqn:weymouth} is referred to as the Weymouth equation for turbulent flow, and its nonconvex nonlinear form presents a critical challenge for optimization applications.

\subsection{Network Modeling}
\label{subsection:network-modeling}
Aside from pipes, gas networks include a variety of other components \citep{koch2015evaluating}, each of which is modeled using different sets of variables and constraints.
This subsection describes the nonconvex models of all components.

\paragraph{\textbf{Notation for Sets.}}
A natural gas transmission network is represented by a directed graph $\mathcal{G} := (\mathcal{N}, \mathcal{A})$, where $\mathcal{N}$ is the set of junctions (nodes) and $\mathcal{A}$ is the set of node-connecting components (arcs).
The set of node-connecting components in the network includes horizontal pipes, short pipes, resistors, loss resistors, valves, regulators, and compressors.
The set of receipts (producers) is denoted by $\mathcal{R}$ and deliveries (consumers) by $\mathcal{D}$.
Receipts and deliveries are attached to existing junctions $i \in \mathcal{N}$.
Furthermore, we let the subset of receipts attached to $i \in \mathcal{N}$ be denoted by $\mathcal{R}_{i}$ and the subset of deliveries by $\mathcal{D}_{i}$.
The sets of horizontal and short pipes in the network are denoted by $\mathcal{P} \subset \mathcal{A}$ and $\mathcal{S} \subset \mathcal{A}$, respectively;
the set of regular and constant loss resistors by $\mathcal{T} \subset \mathcal{A}$ and $\mathcal{U} \subset \mathcal{A}$, respectively;
the set of valves and regulating (i.e., control) valves by $\mathcal{V} \subset \mathcal{A}$ and $\mathcal{W} \subset \mathcal{A}$, respectively;
and the set of compressors by $\mathcal{C} \subset \mathcal{A}$.
Finally, the set of node-connecting components incident to junction $i \in \mathcal{N}$ where $i$ is the tail (respectively, head) of the arc is denoted by $\delta^{+}_{i} := \{(i, j) \in \mathcal{A}\}$ (respectively, $\delta^{-}_{i} := \{(j, i) \in \mathcal{A}\}$).
We now examine each of the aforementioned components individually, define the decision variables, and present constraints that each component enforces on the gas network's operations.
In particular, for each component we present two types of constraints: (i) operational limits and (ii) physical constraints.
We start by examining junctions of the pipeline network.

\paragraph{\textbf{Junctions.}}
Each junction $i \in\mathcal N$ in the network is associated with a pressure variable, $p_i$.
Operational limits require that this pressure resides between predefined lower and upper bounds, denoted by $\underline{p}_{i}$ and $\overline{p}_{i}$, respectively, i.e.,
\begin{equation}
    0 \leq \underline{p}_{i} \leq p_{i} \leq \overline{p}_{i}, ~ \forall i \in \mathcal{N} \label{eqn:pressure-bounds}.
\end{equation}
We note that for a subset of predefined ``slack junctions'' $\mathcal{N}^{s} \subset \mathcal{N}$, $\underline{p}_{i} = \overline{p}_{i}$. 

\paragraph{\textbf{Node-connecting components.}}
Every node-connecting component $(i, j) \in \mathcal A$ is associated with a decision variable, $f_{ij}$, which denotes the mass flow rate across that component. These variables satisfy the capacity constraints
\begin{equation}
    \underline{f}_{ij} \leq f_{ij} \leq \overline{f}_{ij}, ~ \forall (i, j) \in \mathcal{A} \label{eqn:mass-flow-bounds}.
\end{equation}
Here, a positive (respectively, negative) value of $f_{ij}$ implies mass flow that is directed from node $i$ to $j$ (respectively, $j$ to $i$).
Furthermore, for a pipe $(i, j) \in \mathcal{P}$, mass flow $f_{ij}$ and mass flux $\phi_{ij}$ are related by $f_{ij} := \frac{\pi}{4} D_{ij}^{2} \phi_{ij}$.
In the forthcoming paragraphs, we present the constraints required for modeling the operational limits and physics of node-connecting components in the network.

\paragraph{\textbf{Pipes.}}
Pipes transport gas throughout a pipeline network.
We suppose that in steady-state flow, each pipe satisfies the Weymouth Equation \eqref{eqn:weymouth}, which relates the mass flow rate $f_{ij}$ through the cross-sectional area of the pipe to the pressures $p_i$ and $p_j$ at the two end-points of the pipe.
That is,
\begin{equation}
    p_{i}^{2} - p_{j}^{2} = w_{ij} f_{ij} \lvert f_{ij} \rvert, ~ \forall (i, j) \in \mathcal{P} \label{eqn:pipe-weymouth},
\end{equation}
where the mass flow resistance $w_{ij}$ is related to $\gamma_{ij}$ via $w_{ij} := (16 \gamma_{ij}) / (\pi^{2} D_{ij}^{4})$.

\paragraph{\textbf{Short Pipes.}}
Short pipes are components that model resistanceless transport of flow between two junctions.
This is equivalent to treating the length term of a pipe's resistance as negligible.
Short pipes thus ensure the equality of pressures at the two junctions, $i$ and $j$, connected by that component, i.e.,
\begin{equation}
    p_{i} - p_{j} = 0, ~ \forall (i, j) \in \mathcal{S} \label{eqn:short-pipe-pressure}.
\end{equation}

\paragraph{\textbf{Resistors.}}
Aside from pressure losses that arise from the Weymouth equation, a variety of other phenomena can also induce pressure loss.
Examples include turbulence in shaped components, effects of measurement devices, curvature of piping, and partially closed valves.
Resistors serve as surrogate modeling tools for representing these other forms of pressure loss.
Losses across resistors are modeled using the Darcy-Weisbach equation \citep{koch2015evaluating},
\begin{equation}
    p_{i} - p_{j} = \tau_{ij} f_{ij} \lvert f_{ij} \rvert, ~ \forall (i, j) \in \mathcal{T} \label{eqn:resistor-darcy-weisbach}.
\end{equation}
Here, the resistance is defined by $\tau_{ij} := (8 \kappa_{ij})/(\pi^{2} D_{ij}^{4} \rho_s)$, where $\kappa_{ij}$ is the resistor's unitless drag factor, $D_{ij}$ is the resistor's diameter, which may be artificial, and $\rho_s$ is the average standard density of gas throughout the network.
Note that like Constraints \eqref{eqn:pipe-weymouth}, Constraints \eqref{eqn:resistor-darcy-weisbach} are also nonconvex nonlinear.

\paragraph{\textbf{Loss Resistors.}}
Loss resistors serve as an alternate form of the previous resistors, where instead of satisfying the Darcy-Weisbach equation, a fixed pressure loss $\xi_{ij} \geq 0$ is incurred across the component \citep{koch2015evaluating}.
Modeling pressure losses for all loss resistors requires imposing the constraints
\begin{subequations}
\begin{align}
    & f_{ij} (p_{i} - p_{j}) \geq 0, ~ \forall (i, j) \in \mathcal{U} \label{eqn:loss-resistor-direction} \\
    & (p_{i} - p_{j})^{2} = \xi_{ij}^{2}, ~ \forall (i, j) \in \mathcal{U} \label{eqn:loss-resistor-pressure-difference}.
\end{align}
\label{eqn:loss-resistor-pressure}%
\end{subequations}
Here, Constraints \eqref{eqn:loss-resistor-direction} ensure that the mass flow across each loss resistor is in the direction of the pressure loss, and each Constraint \eqref{eqn:loss-resistor-pressure-difference} relates the pressure loss magnitude $\xi_{ij}$ across the loss resistor to the difference of pressures.

\paragraph{\textbf{Valves.}}
Valves are used to route the flow of gas to certain portions of the network or to block flow during maintenance of subnetworks.
In this study, valves are considered controllable elements that are assumed to be either closed or open.
In practice, valves can be partially closed to control the gas velocity, but in these cases, we choose to model the valve as a resistor \citep{koch2015evaluating}.
The operating status of each valve $(i, j) \in \mathcal{V}$ is indicated using a binary variable $z_{ij} \in \{0, 1\}$, where $z_{ij} = 1$ corresponds to an open valve and $z_{ij} = 0$ to a closed valve.
These variables constrain the mass flow across each valve as
\begin{equation}
    \underline{f}_{ij} z_{ij} \leq f_{ij} \leq \overline{f}_{ij} z_{ij}, ~ z_{ij} \in \{0, 1\}, ~ \forall (i, j) \in \mathcal{V} \label{eqn:valve-flow-bounds}.
\end{equation}
Furthermore, when a valve is open, the pressures at the junctions connected by that valve are equal.
When the valve is closed, these pressures are decoupled.
This phenomenon is modeled via the following set of disjunctive constraints:
\begin{subequations}
\begin{align}
    & p_{i} \leq p_{j} + (1 - z_{ij}) \overline{p}_{i}, ~ \forall (i, j) \in \mathcal{V} \label{eqn:valve-pressure-1} \\
    & p_{j} \leq p_{i} + (1 - z_{ij}) \overline{p}_{j}, ~ \forall (i, j) \in \mathcal{V} \label{eqn:valve-pressure-2}.
\end{align}
\label{eqn:valve-pressure}
\end{subequations}

\paragraph{\textbf{Regulators.}}
Large pipes are usually operated at higher pressures than other portions of the network.
As such, interconnection of large pipes with smaller pipes often requires the use of pressure regulators (i.e., control valves) to reduce pressure between differently-sized pipes.
Regulators can also be used as an additional means of controlling flow throughout the network \citep{koch2015evaluating}.
The operating status of a regulator is given using a binary variable $z_{ij} \in \{0, 1\}$, where $z_{ij} = 1$ and $z_{ij} = 0$ indicate active and inactive statuses, respectively.
The mass flow across each regulator is then governed by
\begin{equation}
    \underline{f}_{ij} z_{ij} \leq f_{ij} \leq \overline{f}_{ij} z_{ij}, ~ z_{ij} \in \{0, 1\}, ~ \forall (i, j) \in \mathcal{W} \label{eqn:regulator-flow-bounds}.
\end{equation}
Furthermore, each regulator $(i, j)$ is associated with a multiplicative scaling factor, $\alpha_{ij}$, that defines the relationship between $p_i$ and $p_j$ when the regulator is active, i.e., $\alpha_{ij} p_{i} = p_{j}$.
This factor is constrained by operating limits as $\underline{\alpha}_{ij} = 0 \leq \alpha_{ij} \leq \overline{\alpha}_{ij} = 1$.
As for valves, the pressures at the junctions connected by a regulator are decoupled when the regulator is inactive.
This is modeled by
\begin{subequations}
\begin{align}
    & f_{ij} (p_{i} - p_{j}) \geq 0, ~ \forall (i, j) \in \mathcal{W} \label{eqn:regulator-direction} \\
    & \underline{\alpha}_{ij} p_{i} \leq p_{j} + (1 - z_{ij}) \underline{\alpha}_{ij} \overline{p}_{i}, ~ \forall (i, j) \in \mathcal{W} \label{eqn:regulator-pressure-1} \\
    & p_{j} \leq \overline{\alpha}_{ij} p_{i} + (1 - z_{ij}) \overline{p}_{j}, ~ \forall (i, j) \in \mathcal{W} \label{eqn:regulator-pressure-2}.
\end{align}
\label{eqn:regulator-pressure}%
\end{subequations}
Here, Constraints \eqref{eqn:regulator-direction} ensure that each mass flow is in the same direction as the pressure loss, while Constraints \eqref{eqn:regulator-pressure-1} and \eqref{eqn:regulator-pressure-2} ensure the pressure at node $j$ resides within the scaled bounds of $p_{i}$ when the control valve is open.

\paragraph{\textbf{Compressors.}}
Each compressor $(i, j) \in \mathcal{C}$ boosts the pressure at the downstream junction $j \in \mathcal{N}$ by a variable scalar $\alpha_{ij}$ and is assumed to have negligible length.
For the networks considered in this paper, bidirectional compressors do not exist, although each compressor may or may not allow for \emph{uncompressed} flow in the reverse direction, i.e., from $j$ to $i$.
In this study, these distinct behaviors of compressors are captured by three conditional sets of constraints.

The first set of constraints describes the behavior of compressors where uncompressed reverse flow is prohibited.
These constraints are presented as
\begin{equation}
    \underline{\alpha}_{ij} p_{i} \leq p_{j} \leq \overline{\alpha}_{ij} p_{i}, ~ \forall (i, j) \in \mathcal{C} : \underline{f}_{ij} \geq 0 \label{eqn:compressor-pressures-1},
\end{equation}
where $\underline{\alpha}_{ij}$ and $\overline{\alpha}_{ij}$, $(i, j) \in \mathcal{C}$, are minimum and maximum compression ratios.

The second set of constraints describes the behavior of compressors where reverse flow \emph{is} allowed and the minimum compression ratio is equal to one:
\begin{subequations}
\begin{align}
    & \underline{\alpha}_{ij} p_{i} \leq p_{j} \leq \overline{\alpha}_{ij} p_{i}, ~ \forall (i, j) \in \mathcal{C} : \underline{f}_{ij} < 0 \land \underline{\alpha}_{ij} = 1 \label{eqn:compressor-pressures-2-1} \\
    & f_{ij} (p_{i} - p_{j}) \leq 0, ~ \forall (i, j) \in \mathcal{C} : \underline{f}_{ij} < 0 \land \underline{\alpha}_{ij} = 1 \label{eqn:compressor-pressures-2-3}.
\end{align}
\label{eqn:compressor-pressures-2}%
\end{subequations}
Here, Constraints \eqref{eqn:compressor-pressures-2-3} ensure that, if $f_{ij} < 0$ (i.e., when there exists reverse flow), then $p_{i} = p_{j}$ (i.e., there is no change in pressure across the compressor).

The final set of constraints describes the behavior of compressors where uncompressed reverse flow is allowed and the minimum compression ratio is \emph{not} equal to one.
In this case, the behavior of the compressor must be modeled disjunctively.
To accomplish this, for each compressor, a binary variable $y_{ij} \in \{0, 1\}$ is introduced to denote the direction of flow through the compressor, where $y_{ij} = 1$ implies flow from $i$ to $j$ and $y_{ij} = 0$ implies flow from $j$ to $i$.
The pressures at the junctions that connect each compressor are then modeled as
\begin{subequations}
\begin{align}
    & y_{ij} \in \{0, 1\}, ~ \forall (i, j) \in \mathcal{C} : \underline{f}_{ij} < 0 \land \underline{\alpha}_{ij} \neq 1 \label{eqn:compressor-pressures-3-0} \\
    & p_{j} \leq \overline{\alpha}_{ij} p_{i} + (1 - y_{ij}) \overline{p}_{j}, ~ \forall (i, j) \in \mathcal{C} : \underline{f}_{ij} < 0 \land \underline{\alpha}_{ij} \neq 1 \label{eqn:compressor-pressures-3-1} \\
    & \underline{\alpha}_{ij} p_{i} \leq p_{j} + (1 - y_{ij}) \overline{p}_{i}, ~ \forall (i, j) \in \mathcal{C} : \underline{f}_{ij} < 0 \land \underline{\alpha}_{ij} \neq 1 \label{eqn:compressor-pressures-3-2} \\
    & p_{i} - p_{j} \leq y_{ij} \overline{p}_{i}, ~ \forall (i, j) \in \mathcal{C} : \underline{f}_{ij} < 0 \land \underline{\alpha}_{ij} \neq 1 \label{eqn:compressor-pressures-3-3} \\
    & p_{j} - p_{i} \leq y_{ij} \overline{p}_{j}, ~ \forall (i, j) \in \mathcal{C} : \underline{f}_{ij} < 0 \land \underline{\alpha}_{ij} \neq 1 \label{eqn:compressor-pressures-3-4}.
\end{align}
\label{eqn:compressor-pressures-3}%
\end{subequations}
Here, Constraints \eqref{eqn:compressor-pressures-3-1} and \eqref{eqn:compressor-pressures-3-2} ensure that, when a compressor's flow is positively directed, the pressure at node $j$ is modeled according to the compression ratio bounds.
Constraints \eqref{eqn:compressor-pressures-3-3} and \eqref{eqn:compressor-pressures-3-4} ensure that, when flow is reversed (i.e., $y_{ij} = 0$), the pressures on both sides of the compressor are equal.

\paragraph{\textbf{Receipts, Deliveries, and Mass Conservation.}}
Receipts and deliveries are points in the network that are attached to junctions where gas can be supplied to and withdrawn from the network, respectively.
Each receipt (respectively, delivery) is associated with a nonnegative constant $\overline{s}_{k}$ (respectively, $\overline{d}_{k}$) that denotes the fixed mass supply (respectively, demand) at receipt $k \in \mathcal{R}$ (respectively, delivery $k \in \mathcal{D}$).
Mass conservation throughout the network then requires nodal balance constraints to be enforced at every junction, namely
\begin{equation}
	\sum_{\mathclap{(i, j) \in \delta^{+}_{i}}} f_{ij} - \sum_{\mathclap{(j, i) \in \delta^{-}_{i}}} f_{ji} = \sum_{\mathclap{k \in \mathcal{R}_{i}}} \overline{s}_{k} - \sum_{\mathclap{k \in \mathcal{D}_{i}}} \overline{d}_{k}, \; \forall i \in \mathcal{N} \label{eqn:mass-conservation}.
\end{equation}

\paragraph{\textbf{Feasibility Problem.}}
Given the constraints that model each component in the gas network, the nonconvex program for steady-state feasibility is thus
\begin{equation}\tag{MINCP-F}
\begin{aligned}
	& \textnormal{Pressure and mass flow bounds: Constraints} ~ \eqref{eqn:pressure-bounds}, \eqref{eqn:mass-flow-bounds} \\
    & \textnormal{Pipe dynamics: Constraints} ~ \eqref{eqn:pipe-weymouth} \\
    & \textnormal{Short pipe dynamics: Constraints} ~ \eqref{eqn:short-pipe-pressure} \\
    & \textnormal{Resistor dynamics: Constraints} ~ \eqref{eqn:resistor-darcy-weisbach} \\
    & \textnormal{Loss resistor dynamics: Constraints} ~ \eqref{eqn:loss-resistor-pressure} \\
    & \textnormal{Valve dynamics: Constraints} ~ \eqref{eqn:valve-flow-bounds}, \eqref{eqn:valve-pressure} \\
    & \textnormal{Regulator dynamics: Constraints} ~ \eqref{eqn:regulator-flow-bounds}, \eqref{eqn:regulator-pressure} \\
    & \textnormal{Compressor dynamics: Constraints} ~ \eqref{eqn:compressor-pressures-1}, \eqref{eqn:compressor-pressures-2}, \eqref{eqn:compressor-pressures-3} \\
    & \textnormal{Conservation of mass flow: Constraints} ~ \eqref{eqn:mass-conservation}.
\end{aligned}\label{eqn:nlp}
\end{equation}
The goal of \eqref{eqn:nlp} is to determine if gas can be routed through the pipeline network while satisfying the operational limits and physical constraints imposed by each component of the network.
As formulated, the problem is a mixed-integer nonconvex nonlinear program.
The nonconvexities of the system of equations \eqref{eqn:nlp} arise from three sources: (i) the discreteness of controllable components; (ii) bilinear products of variables appearing in flow direction-related inequalities, i.e., Constraints \eqref{eqn:loss-resistor-direction}, \eqref{eqn:regulator-direction}, and \eqref{eqn:compressor-pressures-2-3}; and (iii) nonlinear equations, i.e., Constraints \eqref{eqn:pipe-weymouth}, \eqref{eqn:resistor-darcy-weisbach}, and \eqref{eqn:loss-resistor-pressure-difference}.
Section \ref{section:mld-formulations} first describes a mixed-integer nonconvex quadratic formulation that addresses (ii), followed by a mixed-integer convex quadratic relaxation that addresses (iii).

\section{Maximal Load Delivery Formulations}
\label{section:mld-formulations}
This section formulates the MLD problem, which seeks to determine a feasible operating point for a damaged gas pipeline network that maximizes the delivery of prioritized load.
To the best of our knowledge, this is the first time such a problem has been formulated for gas networks.
Three variants of the problem are presented in order of successive reformulations and relaxations that, though they increase problem complexity, decrease nonconvexity and hence computational difficulty.
First, Section \ref{subsection:minqp-formulation} extends the nonconvex feasibility constraints of \eqref{eqn:nlp} to formulate the initial MLD problem.
Section \ref{subsection:minqp-reformulation} introduces an \emph{exact} mixed-integer nonconvex quadratic reformulation of this problem.
Finally, to alleviate the challenges associated with nonconvexity, a mixed-integer convex quadratic \emph{relaxation} is presented in Section \ref{subsection:micqp-relaxation}.

Aside from formulating the MLD problem, compared to previous studies (e.g., \citealp{borraz2016convex,wu2017adaptive,chen2018steady,ahumada2021nk,bent2018joint}), this section also (i) presents mixed-integer quadratic reformulations for a broader set of gas pipeline network components and (ii) introduces convex relaxation techniques for components that were not explicitly considered in these studies (e.g., resistors and loss resistors).

\subsection{Mixed-integer Nonconvex Formulation}
\label{subsection:minqp-formulation}
A damaged gas network requires the exclusion of components from the model described in Section \ref{subsection:network-modeling}.
This set of excluded components comprises both the damaged components themselves as well as any connected components.
For example, a damaged junction $i \in \mathcal{N}$ implies a nonoperational status of all node-connecting components $\delta_{i}^{+} \cup \delta_{i}^{-}$.
To this end, we introduce a tilde-based notation to define the sets of components that are \emph{still functional} within a damaged gas network, e.g., $\widetilde{\mathcal{P}} \subseteq \mathcal{P}$ denotes the set of still-operational pipes.

The fundamental motivation for formulating the MLD problem is that a damaged gas network may not be able to satisfy all demands of the original system.
This implies a possible imbalance of the mass conservation Constraints \eqref{eqn:mass-conservation}.
To address this, in the MLD problem, all receipts and deliveries are treated as dispatchable.
This implies the previous constant supplies and demands, $\overline{s}_{k}$, $k \in \mathcal{R}_{i}$, and $\overline{d}_{k}$, $k \in \mathcal{D}_{i}$, for $i \in \mathcal{N}$, become variables bounded as
\begin{subequations}%
\begin{align}%
    0 \leq s_{k} \leq \overline{s}_{k}, ~ &\forall k \in \mathcal{R}_{i}, ~ \forall i \in \widetilde{\mathcal{N}} \\
    0 \leq d_{k} \leq \overline{d}_{k}, ~ &\forall k \in \mathcal{D}_{i}, ~ \forall i \in \widetilde{\mathcal{N}},%
\end{align}%
\label{eqn:dispatchable-bounds}%
\end{subequations}
where $\overline{s}_{k}$ and $\overline{d}_{k}$ denote the original supplies and demands, respectively.
Using these variable loads, the previous mass conservation Constraints \eqref{eqn:mass-conservation} become
\begin{equation}
    \sum_{\mathclap{(i, j) \in \widetilde{\delta}^{+}_{i}}} f_{ij} - \sum_{\mathclap{(j, i) \in \widetilde{\delta}^{-}_{i}}} f_{ji} = \sum_{\mathclap{k \in \mathcal{R}_{i}}} s_{k} - \sum_{\mathclap{k \in \mathcal{D}_{i}}} d_{k}, \; \forall i \in \widetilde{\mathcal{N}} \label{eqn:mass-conservation-dispatch}.
\end{equation}
Next, parameters $\beta_{k} \geq 0$, $k \in \mathcal{D}$, are introduced to denote load restoration priorities.
If no load priorities are available, values of one can be used instead.
(Indeed, this parameterization is used for all of our experiments in Section \ref{section:computational_experiments}.)
Maximization of prioritized load delivered then implies the objective function
\begin{equation}
    \eta(d) = \sum_{i \in \widetilde{\mathcal{N}}} \sum_{k \in \mathcal{D}_{i}} \beta_{k} d_{k} \label{eqn:mld-objective}.
\end{equation}

The mixed-integer nonconvex, nonlinear MLD problem formulation is then
\begin{equation}\tag{MINCP}\begin{aligned}
    & \text{maximize}
	& & \textnormal{Objective function:} ~ \eta(d) ~ \textnormal{of Equation} ~ \eqref{eqn:mld-objective} \\
    & \text{subject to}
    & & \textnormal{Supply and demand bounds: Constraints} ~ \eqref{eqn:dispatchable-bounds} \\
    & & & \textnormal{Conservation of mass flow: Constraints} ~ \eqref{eqn:mass-conservation-dispatch} \\
    & & & \eqref{eqn:nlp} ~ \textnormal{without Constraints} ~ \eqref{eqn:mass-conservation},
\end{aligned}\label{eqn:minlp}\end{equation}
where in \eqref{eqn:nlp}, all component sets are assumed to be replaced with their tilde-denoted counterparts to indicate the application of a damage scenario.

\subsection{Mixed-integer Nonconvex Quadratic Reformulation}
\label{subsection:minqp-reformulation}
As described at the end of Section \ref{subsection:network-modeling}, one source of nonconvexity in \eqref{eqn:minlp} is the existence of bilinear variable products appearing in flow direction-related inequalities.
To address these nonconvexities, this subsection introduces (i) binary direction variables $y_{ij} \in \{0, 1\}$ for all $(i, j) \in \widetilde{\mathcal{A}}$ and (ii) squared pressure variables $\pi_{i} = p_{i}^{2}$ for all $i \in \widetilde{\mathcal{N}}$, which allows for the construction of an \emph{exact} mixed-integer nonconvex quadratic reformulation of the original problem.
That is, the new reformulation contains fewer nonlinearities than \eqref{eqn:minlp} but has the same solution set.
This new formulation enables \emph{global} solutions to be found with modern mixed-integer quadratic programming solvers (namely, \textsc{Gurobi}).
Later, in Section \ref{subsection:micqp-relaxation}, this directed formulation is relaxed to form a mixed-integer \emph{convex} quadratic relaxation of the original MLD problem.

\paragraph{\textbf{Junctions.}}
Squared pressures must reside between predefined bounds, i.e.,
\begin{equation}
    \underline{\pi}_{i} \leq \pi_{i} \leq \overline{\pi}_{i}, ~ \forall i \in \widetilde{\mathcal{N}} \label{eqn:squared-pressure-bounds}.
\end{equation}
where $\underline{\pi}_{i}$ and $\overline{\pi}_{i}$, $i \in \widetilde{\mathcal{N}}$ are derived from the bounds $\underline{p}_{i}$ and $\overline{p}_{i}$, respectively.

\paragraph{\textbf{Node-connecting Components.}}
Using the binary flow direction variables, the mass flow bounds of all node-connecting components are restricted as
\begin{equation}
    (1 - y_{ij}) \underline{f}_{ij} \leq f_{ij} \leq y_{ij} \overline{f}_{ij}, ~ y_{ij} \in \{0, 1\}, ~ \forall (i, j) \in \widetilde{\mathcal{A}} \label{eqn:mass-flow-bounds-miqp}.
\end{equation}

\paragraph{\textbf{Pipes.}}
Squared pressure variables and directions allow for the reduction of nonlinearities in the Weymouth Constraints \eqref{eqn:pipe-weymouth} for pipes, rewritten here as
\begin{subequations}%
\begin{align}%
    \pi_{i} - \pi_{j} &\geq w_{ij} f_{ij}^{2} - (1 - y_{ij}) \left[w_{ij} \underline{f}_{ij}^{2} - (\underline{\pi}_{i} - \overline{\pi}_{j})\right], ~ \forall (i, j) \in \widetilde{\mathcal{P}} \label{eqn:weymouth-pipe-miqp-1} \\
    \pi_{i} - \pi_{j} &\leq w_{ij} f_{ij}^{2}, ~ \forall (i, j) \in \widetilde{\mathcal{P}} \label{eqn:weymouth-pipe-miqp-2} \\
    \pi_{j} - \pi_{i} &\geq w_{ij} f_{ij}^{2} - y_{ij} \left[w_{ij} \overline{f}_{ij}^{2} - (\underline{\pi}_{j} - \overline{\pi}_{i})\right], ~ \forall (i, j) \in \widetilde{\mathcal{P}} \label{eqn:weymouth-pipe-miqp-3} \\
    \pi_{j} - \pi_{i} &\leq w_{ij} f_{ij}^{2}, ~ \forall (i, j) \in \widetilde{\mathcal{P}} \label{eqn:weymouth-pipe-miqp-4}.
\end{align}%
\label{eqn:weymouth-pipe-miqp}%
\end{subequations}
This is similar to the reformulation presented by \citet{borraz2016convex}.
Here, each Constraint \eqref{eqn:weymouth-pipe-miqp-1} ensures a pressure decrease from $i$ to $j$ when $y_{ij} = 1$ (nonnegative flow).
Constraint \eqref{eqn:weymouth-pipe-miqp-2} ensures the Weymouth equation is satisfied when $y_{ij} = 1$.
Constraint \eqref{eqn:weymouth-pipe-miqp-3} ensures a pressure decrease from $j$ to $i$ when $y_{ij} = 0$.
Constraint \eqref{eqn:weymouth-pipe-miqp-4} ensures the Weymouth equation is satisfied when $y_{ij} = 0$.
Directions are used to bound squared pressure differences via
\begin{subequations}%
\begin{align}%
    & (1 - y_{ij}) (\underline{\pi}_{i} - \overline{\pi}_{j}) \leq \pi_{i} - \pi_{j}, ~ \forall (i, j) \in \widetilde{\mathcal{P}} \label{eqn:pipe-direction-miqp-1} \\
    & \pi_{i} - \pi_{j} \leq y_{ij} (\overline{\pi}_{i} - \underline{\pi}_{j}), ~ \forall (i, j) \in \widetilde{\mathcal{P}} \label{eqn:pipe-direction-miqp-2}.
\end{align}%
\label{eqn:pipe-direction-miqp}%
\end{subequations}

\paragraph{\textbf{Short Pipes.}}
Squared pressures allow for rewriting Constraints \eqref{eqn:short-pipe-pressure} as
\begin{equation}
    \pi_{i} - \pi_{j} = 0, ~ \forall (i, j) \in \widetilde{\mathcal{S}} \label{eqn:short-pipe-loss-miqp}.
\end{equation}
As with pipes, direction variables are used to bound squared pressures via
\begin{subequations}%
\begin{align}%
    & (1 - y_{ij}) (\underline{\pi}_{i} - \overline{\pi}_{j}) \leq \pi_{i} - \pi_{j}, ~ \forall (i, j) \in \widetilde{\mathcal{S}} \label{eqn:short-pipe-direction-miqp-1} \\
    & \pi_{i} - \pi_{j} \leq y_{ij} (\overline{\pi}_{i} - \underline{\pi}_{j}), ~ \forall (i, j) \in \widetilde{\mathcal{S}} \label{eqn:short-pipe-direction-miqp-2}.
\end{align}%
\label{eqn:short-pipe-direction-miqp}%
\end{subequations}

\paragraph{\textbf{Resistors.}}
For junctions that are connected to a resistor, the squared pressure variables $\pi_{i}$, $i \in \widetilde{\mathcal{N}}$, must be related to variables denoting pressure, i.e.,
\begin{equation}%
    p_{i}^{2} - \pi_{i} = 0, ~ i \in \widetilde{\mathcal{N}} : \left(\exists (i, j) \in \widetilde{\mathcal{T}}\right) \lor \left(\exists (j, i) \in \widetilde{\mathcal{T}}\right) \label{eqn:resistor-pressure-miqp}.
\end{equation}
The potential loss is then written in a manner as for that of pipes, i.e.,
\begin{subequations}%
\begin{align}%
    p_{i} - p_{j} &\geq w_{ij} f_{ij}^{2} - (1 - y_{ij}) \left[w_{ij} \underline{f}_{ij}^{2} - (\underline{p}_{i} - \overline{p}_{j})\right], ~ \forall (i, j) \in \widetilde{\mathcal{T}} \label{eqn:darcy-weisbach-resistor-miqp-1} \\
    p_{i} - p_{j} &\leq w_{ij} f_{ij}^{2}, ~ \forall (i, j) \in \widetilde{\mathcal{T}} \label{eqn:darcy-weisbach-resistor-miqp-2} \\
    p_{j} - p_{i} &\geq w_{ij} f_{ij}^{2} - y_{ij} \left[w_{ij} \overline{f}_{ij}^{2} - (\underline{p}_{j} - \overline{p}_{i})\right], ~ \forall (i, j) \in \widetilde{\mathcal{T}} \label{eqn:darcy-weisbach-resistor-miqp-3} \\
    p_{j} - p_{i} &\leq w_{ij} f_{ij}^{2}, ~ \forall (i, j) \in \widetilde{\mathcal{T}} \label{eqn:darcy-weisbach-resistor-miqp-4},
\end{align}%
\label{eqn:darcy-weisbach-resistor-miqp}%
\end{subequations}
Direction variables are also used to bound the original pressure variables via
\begin{subequations}%
\begin{align}%
    & (1 - y_{ij}) (\underline{p}_{i} - \overline{p}_{j}) \leq p_{i} - p_{j}, ~ \forall (i, j) \in \widetilde{\mathcal{T}} \label{eqn:resistor-direction-miqp-1} \\
    & p_{i} - p_{j} \leq y_{ij} (\overline{p}_{i} - \underline{p}_{j}), ~ \forall (i, j) \in \widetilde{\mathcal{T}} \label{eqn:resistor-direction-miqp-2}.
\end{align}%
\label{eqn:resistor-direction-miqp}%
\end{subequations}

\paragraph{\textbf{Loss Resistors.}}
As for resistors, squared pressure variables $\pi_{i}$, $i \in \widetilde{\mathcal{N}}$, for junctions connected to a loss resistor, must be related to pressure variables via
\begin{equation}%
    p_{i}^{2} - \pi_{i} = 0, ~ i \in \widetilde{\mathcal{N}} : \left(\exists (i, j) \in \widetilde{\mathcal{U}}\right) \lor \left(\exists (j, i) \in \widetilde{\mathcal{U}}\right) \label{eqn:loss-resistor-pressure-miqp}.
\end{equation}
Note that Constraints \eqref{eqn:loss-resistor-pressure} can be rewritten using absolute values as
\begin{equation}
    \xi_{ij} = \lvert p_{i} - p_{j} \rvert, ~ \forall (i, j) \in \widetilde{\mathcal{U}}.
\end{equation}
However, directions $y_{ij}$ allow this disjunctive form to be modeled linearly, i.e.,
\begin{equation}
    (2 y_{ij} - 1) \xi_{ij} = p_{i} - p_{j}, ~ \forall (i, j) \in \widetilde{\mathcal{U}} \label{eqn:loss-resistor-pressure-sqr-miqp}.
\end{equation}
Here, the direction of pressure loss coincides with the direction given by $y_{ij}$.

\paragraph{\textbf{Valves.}}
Pressure Constraints \eqref{eqn:valve-pressure} are written with squared pressures as
\begin{subequations}
\begin{align}
    & \pi_{i} \leq \pi_{j} + (1 - z_{ij}) \overline{\pi}_{i}, ~ \forall (i, j) \in \widetilde{\mathcal{V}} \label{eqn:valve-pressure-miqp-1} \\
    & \pi_{j} \leq \pi_{i} + (1 - z_{ij}) \overline{\pi}_{j}, ~ \forall (i, j) \in \widetilde{\mathcal{V}} \label{eqn:valve-pressure-miqp-2}.
\end{align}
\label{eqn:valve-pressure-miqp}
\end{subequations}

\paragraph{\textbf{Regulators.}}
Pressure Constraints \eqref{eqn:regulator-pressure} are written with squared pressures as
\begin{subequations}
\begin{align}
    & \pi_{j} - \overline{\alpha}_{ij}^{2} \pi_{i} \leq (2 - y_{ij} - z_{ij}) \overline{\pi}_{j}, ~ \forall (i, j) \in \widetilde{\mathcal{W}} \label{eqn:regulator-pressure-miqp-1} \\
    & \underline{\alpha}_{ij}^{2} \pi_{i} - \pi_{j} \leq (2 - y_{ij} - z_{ij}) \overline{\pi}_{i}, ~ \forall (i, j) \in \widetilde{\mathcal{W}} \label{eqn:regulator-pressure-miqp-2} \\
    & \pi_{j} - \pi_{i} \leq (1 + y_{ij} - z_{ij}) \overline{\pi}_{j}, ~ \forall (i, j) \in \widetilde{\mathcal{W}} \label{eqn:regulator-pressure-miqp-3} \\
    & \pi_{i} - \pi_{j} \leq (1 + y_{ij} - z_{ij}) \overline{\pi}_{i}, ~ \forall (i, j) \in \widetilde{\mathcal{W}} \label{eqn:regulator-pressure-miqp-4}.
\end{align}%
\label{eqn:regulator-pressure-miqp}%
\end{subequations}
Here, when $y_{ij} = 1$ and $z_{ij} = 1$ (i.e., the regulating valve is open and the flow direction is positive), Constraints \eqref{eqn:regulator-pressure-miqp-1} and \eqref{eqn:regulator-pressure-miqp-2} ensure that $\pi_{j}$ resides between the scaled value of $\pi_{i}$.
When $y_{ij} = 0$ and $z_{ij} = 1$, Constraints \eqref{eqn:regulator-pressure-miqp-3} and \eqref{eqn:regulator-pressure-miqp-4} ensure that $\pi_{i} = \pi_{j}$, and reverse flow is allowed.
Finally, when the valve is closed, $z_{ij} = 0$ and $y_{ij} \in \{0, 1\}$, which ensures $\pi_{i}$ and $\pi_{j}$ are decoupled.

\paragraph{\textbf{Compressors.}}
Constraints \eqref{eqn:compressor-pressures-1}, which model compressors where uncompressed reverse flow is prohibited, are written with squared pressures as
\begin{equation}
    \underline{\alpha}_{ij}^{2} \pi_{i} \leq \pi_{j} \leq \overline{\alpha}_{ij}^{2} \pi_{i}, ~ \forall (i, j) \in \mathcal{C} : \underline{f}_{ij} \geq 0 \label{eqn:compressor-pressures-miqp-1}.
\end{equation}
For compressors that allow reverse flow, Constraints \eqref{eqn:compressor-pressures-2} and \eqref{eqn:compressor-pressures-3} become
\begin{subequations}
\begin{align}
    & \pi_{j} \leq \overline{\alpha}_{ij}^{2} \pi_{i} + (1 - y_{ij}) \overline{\pi}_{j}, ~ \forall (i, j) \in \mathcal{C} : \underline{f}_{ij} < 0 \label{eqn:compressor-pressures-miqp-2-1} \\
    & \underline{\alpha}_{ij}^{2} \pi_{i} \leq \pi_{j} + (1 - y_{ij}) \underline{\alpha}_{ij}^{2} \overline{\pi}_{i}, ~ \forall (i, j) \in \mathcal{C} : \underline{f}_{ij} < 0 \label{eqn:compressor-pressures-miqp-2-2} \\
    & \pi_{i} - \pi_{j} \leq y_{ij} \overline{\pi}_{i}, ~ \forall (i, j) \in \mathcal{C} : \underline{f}_{ij} < 0 \label{eqn:compressor-pressures-miqp-2-3} \\
    & \pi_{j} - \pi_{i} \leq y_{ij} \overline{\pi}_{j}, ~ \forall (i, j) \in \mathcal{C} : \underline{f}_{ij} < 0 \label{eqn:compressor-pressures-miqp-2-4}.
\end{align}
\label{eqn:compressor-pressures-miqp-2}%
\end{subequations}
Here, when $y_{ij} = 1$, Constraints \eqref{eqn:compressor-pressures-miqp-2-1} and \eqref{eqn:compressor-pressures-miqp-2-2} require $\pi_{j}$ to reside within the scaled bounds of $\pi_{i}$.
When $y_{ij} = 0$, Constraints \eqref{eqn:compressor-pressures-miqp-2-3} and \eqref{eqn:compressor-pressures-miqp-2-4} ensure the equality of pressures when flow is from $j$ to $i$ (i.e., there is no compression).

\paragraph{\textbf{Direction-related Valid Inequalities.}}
The formulations of \citet{borraz2016convex} include valid inequalities that relate node-connecting component directions $y_{ij}$ to nodal conditions in the directed network.
These inequalities improve relaxed formulations.
Here, we employ the constraints that model flow directionality at junctions with zero supply, zero demand, and degree two:
\begin{subequations}\begin{align}
    & \sum_{\mathclap{(i, j) \in \widetilde{\delta}_{i}^{-}}} y_{ij} - \sum_{\mathclap{(i, j) \in \widetilde{\delta}_{i}^{+}}} y_{ij} = 0, ~ i \in \widetilde{\mathcal{N}}_{0} : \left(\left\lvert \widetilde{\delta}_{i}^{\pm} \right\rvert = 1\right) \label{eqn:deg-2-flow-1} \\
    & \sum_{\mathclap{(i, j) \in \widetilde{\delta}_{i}^{-}}} y_{ij} + \sum_{\mathclap{(i, j) \in \widetilde{\delta}_{i}^{+}}} y_{ij} = 1, ~ i \in \widetilde{\mathcal{N}}_{0} : \left(\left\lvert\widetilde{\delta}_{i}^{\pm}\right\rvert = 2\right) \land \left(\left\lvert \widetilde{\delta}_{i}^{\mp} \right\rvert = 0\right) \label{eqn:deg-2-flow-2},
\end{align}\label{eqn:deg-2-flow}\end{subequations}
where $\widetilde{\mathcal{N}}_{0} \subset \widetilde{\mathcal{N}}$ denotes the subset of junctions where $\sum_{k \in \mathcal{R}_{i}} \overline{s}_{i} = \sum_{k \in \mathcal{D}_{i}} \overline{d}_{i} = 0$.
These constraints imply that, for this subset of junctions, the direction of incoming mass flow must be equal to the direction of outgoing mass flow.

\paragraph{\textbf{Reformulation.}}
The preceding changes enable the mixed-integer nonconvex \eqref{eqn:minlp} to be formulated as the mixed-integer \emph{nonconvex quadratic} program
\begin{equation}\tag{MINQP}\begin{aligned}
    & \text{maximize}
	& & \textnormal{Objective function:} ~ \eta(d) ~ \textnormal{of Equation} ~ \eqref{eqn:mld-objective} \\
    & \text{subject to}
    & & \textnormal{Supply and demand bounds: Constraints} ~ \eqref{eqn:dispatchable-bounds} \\
    & & & \textnormal{Conservation of mass flow: Constraints} ~ \eqref{eqn:mass-conservation-dispatch} \\
    & & & \textnormal{Pressure bounds: Constraints} ~ \eqref{eqn:pressure-bounds}, \eqref{eqn:squared-pressure-bounds} \\
    & & & \textnormal{Directed mass flow bounds: Constraints} ~ \eqref{eqn:mass-flow-bounds-miqp} \\
    & & & \textnormal{Pipe dynamics: Constraints} ~ \eqref{eqn:weymouth-pipe-miqp}, \eqref{eqn:pipe-direction-miqp} \\
    & & & \textnormal{Short pipe dynamics: Constraints} ~ \eqref{eqn:short-pipe-loss-miqp}, \eqref{eqn:short-pipe-direction-miqp} \\
    & & & \textnormal{Resistor dynamics: Constraints} ~ \eqref{eqn:resistor-pressure-miqp}, \eqref{eqn:darcy-weisbach-resistor-miqp}, \eqref{eqn:resistor-direction-miqp} \\
    & & & \textnormal{Loss resistor dynamics: Constraints} ~ \eqref{eqn:loss-resistor-pressure-miqp}, \eqref{eqn:loss-resistor-pressure-sqr-miqp} \\
    & & & \textnormal{Valve dynamics: Constraints} ~ \eqref{eqn:valve-flow-bounds}, \eqref{eqn:valve-pressure-miqp} \\
    & & & \textnormal{Regulator dynamics: Constraints} ~ \eqref{eqn:regulator-flow-bounds}, \eqref{eqn:regulator-pressure-miqp} \\
    & & & \textnormal{Compressor dynamics: Constraints} ~ \eqref{eqn:compressor-pressures-miqp-1}, \eqref{eqn:compressor-pressures-miqp-2} \\
    & & & \textnormal{Direction-related cuts: Constraints} ~ \eqref{eqn:deg-2-flow}.
\end{aligned}\label{eqn:minqp}\end{equation}

\subsection{Mixed-integer Convex Quadratic Relaxation}
\label{subsection:micqp-relaxation}
There are a significant number of nonconvex nonlinear constraints in the problem \eqref{eqn:minqp} that cause it to become computationally intractable with increasing network size.
Specifically, these are Constraints \eqref{eqn:weymouth-pipe-miqp-2}, \eqref{eqn:weymouth-pipe-miqp-4}, \eqref{eqn:resistor-pressure-miqp}, \eqref{eqn:darcy-weisbach-resistor-miqp-2}, \eqref{eqn:darcy-weisbach-resistor-miqp-4}, and \eqref{eqn:loss-resistor-pressure-miqp}.
In this section, we apply convex relaxation strategies, as done in \citet{borraz2016convex}, \citet{wu2017adaptive}, and \citet{chen2018steady}, to address this.
We also extend these studies by formulating relaxations for components that were not previously considered, e.g., resistors and loss resistors, which require constraints involving explicit pressure variables.

\paragraph{\textbf{Pipes.}}
We first apply convex relaxations to the Weymouth Constraints \eqref{eqn:weymouth-pipe-miqp} for pipes.
The variables $\ell_{ij}$ for $(i, j) \in \widetilde{\mathcal{P}}$ are introduced to denote the difference in squared pressures across each pipe.
The relaxation constraints are then
\begin{subequations}%
\begin{align}%
    & \pi_{j} - \pi_{i} \leq \ell_{ij} \leq \pi_{i} - \pi_{j}, ~ \forall (i, j) \in \widetilde{\mathcal{P}} \label{eqn:weymouth-pipe-micqp-1} \\
    & \ell_{ij} \leq \pi_{j} - \pi_{i} + (2 y_{ij}) (\overline{\pi}_{i} - \underline{\pi}_{j}), ~ \forall (i, j) \in \widetilde{\mathcal{P}} \label{eqn:weymouth-pipe-micqp-2} \\
    & \ell_{ij} \leq \pi_{i} - \pi_{j} + (2 y_{ij} - 2) (\underline{\pi}_{i} - \overline{\pi}_{j}), ~ \forall (i, j) \in \widetilde{\mathcal{P}} \label{eqn:weymouth-pipe-micqp-3} \\
    & w_{ij} f_{ij}^{2} \leq \ell_{ij}, ~ \forall (i, j) \in \widetilde{\mathcal{P}} \label{eqn:weymouth-pipe-micqp-4} \\
    & \ell_{ij} \leq w_{ij} \overline{f}_{ij} f_{ij} + (1 - y_{ij}) \left(\left\lvert w_{ij} \overline{f}_{ij} \underline{f}_{ij}\right\rvert + w_{ij} \underline{f}_{ij}^{2}\right), ~ \forall (i, j) \in \widetilde{\mathcal{P}} \label{eqn:weymouth-pipe-micqp-5} \\
    & \ell_{ij} \leq w_{ij} \underline{f}_{ij} f_{ij} + y_{ij} \left(\left\lvert w_{ij} \overline{f}_{ij} \underline{f}_{ij}\right\rvert + w_{ij} \overline{f}_{ij}^{2}\right), ~ \forall (i, j) \in \widetilde{\mathcal{P}} \label{eqn:weymouth-pipe-micqp-6}.
\end{align}%
\label{eqn:weymouth-pipe-micqp}%
\end{subequations}
Here, Constraints \eqref{eqn:weymouth-pipe-micqp-1} ensure each loss of squared pressures resides between the corresponding differences.
Constraints \eqref{eqn:weymouth-pipe-micqp-2} ensure that when $y_{ij} = 0$, each loss is bounded by $\pi_{j} - \pi_{i}$, and
Constraints \eqref{eqn:weymouth-pipe-micqp-3} imply that each loss is bounded by $\pi_{i} - \pi_{j}$ when $y_{ij} = 1$.
Constraints \eqref{eqn:weymouth-pipe-micqp-4} are the primarily convex relaxations of the Weymouth equations.
Finally, Constraints \eqref{eqn:weymouth-pipe-micqp-5} apply linear upper bounds on each variable $\ell_{ij}$, depending on the choice of flow direction.

\paragraph{\textbf{Resistors.}}
We next apply convex relaxations to Constraints \eqref{eqn:resistor-pressure-miqp}, which relate nonsquared and squared pressure variables.
These relaxations yield
\begin{equation}%
    p_{i}^{2} \leq \pi_{i}, ~ i \in \widetilde{\mathcal{N}} : \left(\exists (i, j) \in \widetilde{\mathcal{T}}\right) \lor \left(\exists (j, i) \in \widetilde{\mathcal{T}}\right) \label{eqn:resistor-pressure-micqp}.
\end{equation}
Constraints \eqref{eqn:darcy-weisbach-resistor-miqp} are then relaxed as done for Constraints \eqref{eqn:weymouth-pipe-micqp}, i.e.,
\begin{subequations}%
\begin{align}%
    & p_{j} - p_{i} \leq \ell_{ij} \leq p_{i} - p_{j}, ~ \forall (i, j) \in \widetilde{\mathcal{T}} \label{eqn:darcy-weisbach-resistor-micqp-1} \\
    & \ell_{ij} \leq p_{j} - p_{i} + (2 y_{ij}) (\overline{p}_{i} - \underline{p}_{j}), ~ \forall (i, j) \in \widetilde{\mathcal{T}} \label{eqn:darcy-weisbach-resistor-micqp-2} \\
    & \ell_{ij} \leq p_{i} - p_{j} + (2 y_{ij} - 2) (\underline{p}_{i} - \overline{p}_{j}), ~ \forall (i, j) \in \widetilde{\mathcal{T}} \label{eqn:darcy-weisbach-resistor-micqp-3} \\
    & w_{ij} f_{ij}^{2} \leq \ell_{ij}, ~ \forall (i, j) \in \widetilde{\mathcal{T}} \label{eqn:darcy-weisbach-resistor-micqp-4} \\
    & \ell_{ij} \leq w_{ij} \overline{f}_{ij} f_{ij} + (1 - y_{ij}) \left(\left\lvert w_{ij} \overline{f}_{ij} \underline{f}_{ij}\right\rvert + w_{ij} \underline{f}_{ij}^{2}\right), ~ \forall (i, j) \in \widetilde{\mathcal{T}} \label{eqn:darcy-weisbach-resistor-micqp-5} \\
    & \ell_{ij} \leq w_{ij} \underline{f}_{ij} f_{ij} + y_{ij} \left(\left\lvert w_{ij} \overline{f}_{ij} \underline{f}_{ij}\right\rvert + w_{ij} \overline{f}_{ij}^{2}\right), ~ \forall (i, j) \in \widetilde{\mathcal{T}} \label{eqn:darcy-weisbach-resistor-micqp-6}.
\end{align}%
\label{eqn:darcy-weisbach-resistor-micqp}%
\end{subequations}

\paragraph{\textbf{Loss Resistors.}}
As we have done for resistors, we apply convex relaxations to Constraints \eqref{eqn:loss-resistor-pressure-miqp} relating nonsquared and squared pressures.
We obtain
\begin{equation}%
    p_{i}^{2} \leq \pi_{i}, ~ i \in \widetilde{\mathcal{N}} : \left(\exists (i, j) \in \widetilde{\mathcal{U}}\right) \lor \left(\exists (j, i) \in \widetilde{\mathcal{U}}\right) \label{eqn:loss-resistor-pressure-micqp}.
\end{equation}

\paragraph{\textbf{Relaxation.}}
The mixed-integer convex quadratic relaxation of \eqref{eqn:minqp} is
\begin{equation}\tag{MICQP}\begin{aligned}
    & \text{maximize}
	& & \textnormal{Objective function:} ~ \eta(d) ~ \textnormal{of Equation} ~ \eqref{eqn:mld-objective} \\
    & \text{subject to}
    & & \textnormal{Supply and demand bounds: Constraints} ~ \eqref{eqn:dispatchable-bounds} \\
    & & & \textnormal{Conservation of mass flow: Constraints} ~ \eqref{eqn:mass-conservation-dispatch} \\
    & & & \textnormal{Pressure bounds: Constraints} ~ \eqref{eqn:pressure-bounds}, \eqref{eqn:squared-pressure-bounds} \\
    & & & \textnormal{Directed mass flow bounds: Constraints} ~ \eqref{eqn:mass-flow-bounds-miqp} \\
    & & & \textnormal{Pipe dynamics: Constraints} ~ \eqref{eqn:pipe-direction-miqp}, \eqref{eqn:weymouth-pipe-micqp} \\
    & & & \textnormal{Short pipe dynamics: Constraints} ~ \eqref{eqn:short-pipe-loss-miqp}, \eqref{eqn:short-pipe-direction-miqp} \\
    & & & \textnormal{Resistor dynamics: Constraints} ~ \eqref{eqn:resistor-direction-miqp}, \eqref{eqn:resistor-pressure-micqp}, \eqref{eqn:darcy-weisbach-resistor-micqp} \\
    & & & \textnormal{Loss resistor dynamics: Constraints} ~ \eqref{eqn:loss-resistor-pressure-sqr-miqp}, \eqref{eqn:loss-resistor-pressure-micqp} \\
    & & & \textnormal{Valve dynamics: Constraints} ~ \eqref{eqn:valve-flow-bounds}, \eqref{eqn:valve-pressure-miqp} \\
    & & & \textnormal{Regulator dynamics: Constraints} ~ \eqref{eqn:regulator-flow-bounds}, \eqref{eqn:regulator-pressure-miqp} \\
    & & & \textnormal{Compressor dynamics: Constraints} ~ \eqref{eqn:compressor-pressures-miqp-1}, \eqref{eqn:compressor-pressures-miqp-2} \\
    & & & \textnormal{Direction-related cuts: Constraints} ~ \eqref{eqn:deg-2-flow}.
\end{aligned}\label{eqn:micqp}\end{equation}
In Section \ref{section:computational_experiments}, we present the results of a comprehensive computational study performed to compare the tractability and application of the exact and relaxed MLD problem formulations, \eqref{eqn:minqp} and \eqref{eqn:micqp}, respectively.
These computations display the practical benefits of employing the relaxation, \eqref{eqn:micqp}.

\section{Computational Experiments}
\label{section:computational_experiments}
This section experimentally analyzes the applicability and computational performance of the \eqref{eqn:minqp} and \eqref{eqn:micqp} MLD formulations presented in Section \ref{section:mld-formulations}.
This informs of us the practical, analytical, and computational trade-offs associated with using the exact and relaxed MLD problem formulations, respectively.
To accomplish this, we consider three different types of damage scenarios: (i) $N{-}1$ or single contingency scenarios, (ii) $N{-}k$ or multi-contingency scenarios, and (iii) earthquake damage scenarios.
Each scenario is intended to simulate common sources of damage to a gas pipeline network.
Although these damage models arise from physically reasonable assumptions, we do not claim to quantify the robustness of the specific networks considered.
Rather, these models serve as proofs of concept for studying three aspects of the MLD problem: (i) understanding the tractability of MLD formulations, (ii) highlighting the qualitative insights given by an MLD analysis, and (iii) providing guidelines for future applications of this work to real-world scenarios.

Both MLD formulations are implemented in the \textsc{Julia} programming language using the mathematical modeling layer \textsc{JuMP}, version 0.21 \citep{dunning2017jump}, and version 0.8 of \textsc{GasModels}, an open-source \textsc{Julia} package for steady-state and transient natural gas network optimization \citep{gasmodels}.
Section \ref{subsection:experimental_setup} describes the instances, computational resources, and parameters used throughout these experiments; Section \ref{subsection:n-1_experiments} compares the efficacy of MLD formulations on $N{-}1$ contingency scenarios for each network; Section \ref{subsection:n-k_experiments} does the same for randomized $N{-}k$ multi-contingency scenarios, where $k$ corresponds to $15\%$ of node-connecting components in each network; and Section \ref{subsection:computational-performance} compares the runtime performance of the formulations over both the $N{-}1$ and $N{-}k$ experiment sets.
Finally, Section \ref{subsection:earthquake_experiments} provides a proof of concept application of the MLD problem to damage scenarios precipitated by deterministic or stochastic natural hazards (in this case, earthquakes).

\subsection{Experimental Test Data \& Setup}
\label{subsection:experimental_setup}
The numerical experiments consider networks of various sizes that appear in literature of natural gas transmission network modeling or are derivable from open data.
These instances are summarized in Table \ref{table:networks}.
Here, the \texttt{Belgium-20} network is derived from the application of \citet{de2000gas}; the North American 154-junction network (i.e., \texttt{NA-154}) is derived by subject matter experts using public data; and \texttt{GasLib} networks are obtained directly from \citet{SABHJKKIOSSS17}.
For the \texttt{GasLib-582} and \texttt{GasLib-4197} networks, the \texttt{nomination\_freezing\_1} and \texttt{nomination\_mild\_0006} delivery and receipt nominations are used, respectively.
Generally, steady-state optimization problems involving most networks can be solved to optimality given a small amount of time (e.g., seconds to minutes).
The notable exception is the \texttt{GasLib-4197} network, the largest system, which requires hours to solve many instances.


\begin{table}[t]
    \begin{center}
        \begin{tabular}{|c|c|c|c|c|c|c|c|c|}
            \hline
            Network & $\lvert \mathcal{N} \rvert$ & $\lvert \mathcal{P} \rvert$ & $\lvert \mathcal{S} \rvert$ & $\lvert \mathcal{T} \rvert$ & $\lvert \mathcal{U} \rvert$ & $\lvert \mathcal{V} \rvert$ & $\lvert \mathcal{W} \rvert$ & $\lvert \mathcal{C} \rvert$ \\ \hline
            \texttt{GasLib-11} & 11 & 8 & 0 & 0 & 0 & 1 & 0 & 2 \\ \hline
            \texttt{Belgian-20} & 20 & 24 & 0 & 0 & 0 & 0 & 0 & 3 \\ \hline
            \texttt{GasLib-24} & 24 & 19 & 1 & 1 & 0 & 0 & 1 & 3 \\ \hline
            \texttt{GasLib-40} & 40 & 39 & 0 & 0 & 0 & 0 & 0 & 6 \\ \hline
            \texttt{GasLib-134} & 134 & 86 & 45 & 0 & 0 & 0 & 1 & 1 \\ \hline
            \texttt{GasLib-135} & 135 & 141 & 0 & 0 & 0 & 0 & 0 & 29 \\ \hline
            \texttt{NA-154} & 154 & 140 & 0 & 0 & 0 & 0 & 0 & 12 \\ \hline
            \texttt{GasLib-582} & 582 & 278 & 269 & 8 & 0 & 26 & 23 & 5 \\ \hline
            \texttt{GasLib-4197} & 4197 & 3537 & 343 & 22 & 6 & 426 & 120 & 12 \\ \hline
        \end{tabular}
    \end{center}
    \caption{Summary of natural gas transmission networks derived from open data.}
    \label{table:networks}
\end{table}

Experiments are performed on the Darwin experimental computing cluster at Los Alamos National Laboratory.
Each optimization computation is provided a wall-clock time of one hour on a node containing two Intel Xeon E5-2695 v4 processors, each with 18 cores @2.10 GHz, and 125 GB of memory.
For solutions of the nonconvex and convex MIQPs, \textsc{Gurobi} 9.0 is used with the parameter \texttt{MIPGap=0.0}.
For experiments employing \eqref{eqn:minqp}, the setting \texttt{NonConvex=2} is used, which allows for global optimization of this formulation.

\subsection{Single Contingency Damage Scenarios}
\label{subsection:n-1_experiments}
The first damage model considered is the single contingency or $N{-}1$ model, where $N$ corresponds to the original number of network components and $N{-}1$ indicates that an individual component is removed (or damaged).
This model can be thought of as a method for simulating the effects of an unscheduled component outage.
Here, it is assumed that both nodal components (i.e., junctions) and node-connecting components (e.g., pipes) can comprise an $N{-}1$ damage scenario.
In our study, this damage model is intended to demonstrate feasibility of the MLD problem for a broad variety of network structures and to validate our network modeling assumptions for damaged gas networks.

\begin{table}[t]
    \begin{center}
    \begin{tabular}{c|c|c|c|c|c|c|}
        \cline{2-7} & \multicolumn{3}{c|}{MINQP} & \multicolumn{3}{c|}{MICQP} \\
        \cline{1-7} \multicolumn{1}{|c|}{Network} & Opt. (\%) & Lim. (\%) & Inf. (\%) & Opt. (\%) & Lim. (\%) & Inf. (\%) \\
        \cline{1-7} \multicolumn{1}{|c|}{\texttt{GasLib-11}} & $100.0$ & $0.0$ & $0.0$ & $100.0$ & $0.0$ & $0.0$ \\
        \cline{1-7} \multicolumn{1}{|c|}{\texttt{Belgium-20}} & $100.0$ & $0.0$ & $0.0$ & $100.0$ & $0.0$ & $0.0$ \\
        \cline{1-7} \multicolumn{1}{|c|}{\texttt{GasLib-24}} & $100.0$ & $0.0$ & $0.0$ & $100.0$ & $0.0$ & $0.0$ \\
        \cline{1-7} \multicolumn{1}{|c|}{\texttt{GasLib-40}} & $100.0$ & $0.0$ & $0.0$ & $100.0$ & $0.0$ & $0.0$ \\
        \cline{1-7} \multicolumn{1}{|c|}{\texttt{GasLib-134}} & $100.0$ & $0.0$ & $0.0$ & $100.0$ & $0.0$ & $0.0$ \\
        \cline{1-7} \multicolumn{1}{|c|}{\texttt{GasLib-135}} & $90.2$ & $9.8$ & $0.0$ & $100.0$ & $0.0$ & $0.0$ \\
        \cline{1-7} \multicolumn{1}{|c|}{\texttt{NA-154}} & $100.0$ & $0.0$ & $0.0$ & $100.0$ & $0.0$ & $0.0$ \\
        \cline{1-7} \multicolumn{1}{|c|}{\texttt{GasLib-582}} & $79.6$ & $20.4$ & $0.0$ & $100.0$ & $0.0$ & $0.0$ \\
        \cline{1-7} \multicolumn{1}{|c|}{\texttt{GasLib-4197}} & $3.2$ & $90.3$ & $6.5$ & $43.8$ & $56.2$ & $0.0$ \\
        \cline{1-7}
    \end{tabular}
    \end{center}
    \caption{Comparison of solver termination statuses over all $N{-}1$ contingency scenarios. Here, ``Opt.'' corresponds to instances where optimality is proven, ``Lim.'' to instances where the time limit is reached, and ``Inf.'' to instances that are claimed to be infeasible.}
    \label{table:n-1}
\end{table}

For each network, the set of all such possible $N{-}1$ scenarios is considered, and the corresponding instances are solved using both the \eqref{eqn:minqp} and \eqref{eqn:micqp} MLD formulations.
The exception is \texttt{GasLib-4197}, which is limited to $1040$ unique scenarios because of cluster restrictions.
Table \ref{table:n-1} displays statistics of solver termination statuses across all scenarios for each network and formulation.
For all networks except \texttt{GasLib-135}, \texttt{GasLib-582}, and \texttt{GasLib-4197}, optimal solutions to all instances are found within the prescribed one hour time limit.
For \texttt{GasLib-135}, the \eqref{eqn:minqp} formulation is unable to prove global optimality on $98$ instances.
Notably, \texttt{GasLib-135} contains the largest number of compressors compared to other networks, and these additional degrees of freedom are the source of computational complexity.
For \texttt{GasLib-582} and \texttt{GasLib-4197}, the large number of \eqref{eqn:minqp} instances that cannot be solved to optimality is likely because of the networks' large sizes.

Comparing the two formulations shows the benefit of using the relaxation-based \eqref{eqn:micqp} approach, which is capable of solving much larger proportions of challenging \texttt{GasLib-135}, \texttt{GasLib-582}, and \texttt{GasLib-4197} instances.
This suggests that \eqref{eqn:micqp}, when compared to the mixed-integer nonconvex \eqref{eqn:minqp}, is often a better candidate for numerically intensive applications.
We also note that, when using the default \textsc{Gurobi} parameters described in Section \ref{subsection:experimental_setup}, $292$ of $1040$ \texttt{GasLib-4197} instances of the \eqref{eqn:minqp} MLD formulation are claimed to be infeasible.
Using the parameter \texttt{NumericFocus=3} for this subset of $292$ instances, however, decreases this number to $68$.
These claimed infeasibilities are likely related to the numerical properties of the \texttt{GasLib-4197} data set rather than our MLD formulations.
Additional future work to preprocess data and rescale constraints is warranted to address these issues.

\subsection{Multi-contingency Damage Scenarios}
\label{subsection:n-k_experiments}
The second damage model is the multi-contingency or $N{-}k$ model, where $k$ corresponds to the number of components that are simultaneously removed from the network.
These scenarios are intended to capture the effects of multimodal network failures.
Here, we consider the removal of only node-connecting components within the generated $N{-}k$ scenarios.
In each scenario, a uniformly random selection of $15\%$ node-connecting components are assumed to be nonoperational.
Heuristically, this proportion of components seems to generate challenging scenarios while providing different maximal load distributions across the networks considered.
For each network, one thousand such scenarios are generated.
If solved, the maximal proportional load delivered in each experiment is then computed as the ratio of the optimal nonprioritized objective in Equation \eqref{eqn:mld-objective} and the maximal load of the undamaged network.


\begin{table}[t]
    \begin{center}
    \begin{tabular}{c|c|c|c|c|c|c|}
        \cline{2-7} & \multicolumn{3}{c|}{MINQP} & \multicolumn{3}{c|}{MICQP} \\
        \cline{1-7} \multicolumn{1}{|c|}{Network} & Opt. (\%) & Lim. (\%) & Inf. (\%) & Opt. (\%) & Lim. (\%) & Inf. (\%) \\
        \cline{1-7} \multicolumn{1}{|c|}{\texttt{GasLib-11}} & $100.0$ & $0.0$ & $0.0$ & $100.0$ & $0.0$ & $0.0$ \\
        \cline{1-7} \multicolumn{1}{|c|}{\texttt{Belgium-20}} & $100.0$ & $0.0$ & $0.0$ & $100.0$ & $0.0$ & $0.0$ \\
        \cline{1-7} \multicolumn{1}{|c|}{\texttt{GasLib-24}} & $100.0$ & $0.0$ & $0.0$ & $100.0$ & $0.0$ & $0.0$ \\
        \cline{1-7} \multicolumn{1}{|c|}{\texttt{GasLib-40}} & $100.0$ & $0.0$ & $0.0$ & $100.0$ & $0.0$ & $0.0$ \\
        \cline{1-7} \multicolumn{1}{|c|}{\texttt{GasLib-134}} & $100.0$ & $0.0$ & $0.0$ & $100.0$ & $0.0$ & $0.0$ \\
        \cline{1-7} \multicolumn{1}{|c|}{\texttt{GasLib-135}} & $99.6$ & $0.4$ & $0.0$ & $100.0$ & $0.0$ & $0.0$ \\
        \cline{1-7} \multicolumn{1}{|c|}{\texttt{NA-154}} & $100.0$ & $0.0$ & $0.0$ & $100.0$ & $0.0$ & $0.0$ \\
        \cline{1-7} \multicolumn{1}{|c|}{\texttt{GasLib-582}} & $99.9$ & $0.1$ & $0.0$ & $100.0$ & $0.0$ & $0.0$ \\
        \cline{1-7} \multicolumn{1}{|c|}{\texttt{GasLib-4197}} & $92.7$ & $7.3$ & $0.0$ & $99.2$ & $0.8$ & $0.0$ \\
        \cline{1-7}
    \end{tabular}
    \end{center}
    \caption{Comparison of termination statuses over $1000$ randomized $N{-}k$ multi-contingency scenarios, where $k$ corresponds to $15\%$ of all edge-type components in the network.}
    \label{table:n-k-15}
\end{table}

Table \ref{table:n-k-15} displays statistics of solver termination statuses across all damage scenarios for each network and formulation.
For all except \texttt{GasLib-135}, \texttt{GasLib-582}, and \texttt{GasLib-4197}, globally optimal solutions to MLD instances are found within the one hour time limit.
For \texttt{GasLib-135} and \texttt{GasLib-582}, the \eqref{eqn:minqp} formulation is unable to prove optimality on four instances and one instance, respectively.
For \texttt{GasLib-4197}, a larger proportion cannot be solved.
Comparing the two formulations again shows the benefit of the relaxation-based approach, which solves larger proportions of challenging instances.
We note that one \texttt{GasLib-4197} \eqref{eqn:minqp} instance is claimed to be infeasible using default \textsc{Gurobi} parameters but is resolved when using \texttt{NumericFocus=3}.

\begin{figure}

\centering
\begin{subfigure}[b]{0.345\textwidth}
    \begin{tikzpicture}
        \begin{axis}[height=4.0cm,ybar,enlargelimits=false,bar width=2.00pt,
                     ylabel=Scenarios Solved ($\%$),xticklabels={,,},xtick={0,20,40,60,80,100},
                     xmajorgrids={false},xmin=-5.0,xmax=115.0,ymin=0.0,
                     legend style={at={(0.97, 0.80)},font=\tiny},anchor=north east,
                     tick pos=left,ymax=55.0]
            \node[fill=white,anchor=north east] at (rel axis cs:0.98,0.97) {\texttt{GasLib-11}};
            \addplot[fill] table [x, y, col sep=comma] {data/hist-COLLECTION.n-k-15-FORMULATION.DWP-NETWORK.GasLib-11.csv};
            \node[fill=white,anchor=north east] at (rel axis cs:0.98,0.97) {\texttt{GasLib-11}};
            \addplot[fill=white] table [x, y, col sep=comma] {data/hist-COLLECTION.n-k-15-FORMULATION.CRDWP-NETWORK.GasLib-11.csv};
            \addlegendentry{\eqref{eqn:minqp}};
            \addlegendentry{\eqref{eqn:micqp}};
        \end{axis}
    \end{tikzpicture}
\end{subfigure}
\begin{subfigure}[b]{0.32\textwidth}
    \hspace{0.5em}
    \begin{tikzpicture}
        \begin{axis}[height=4.0cm,ybar,enlargelimits=false,bar width=2.00pt,
                     xmajorgrids={false},xmin=-5.0,xmax=115.0,ymin=0.0,
                     ymax=55.0,xticklabels={,,},xtick={0,20,40,60,80,100},
                     tick pos=left,yticklabels={,,}]
            \node[fill=white,anchor=north east] at (rel axis cs:0.98,0.97) {\texttt{Belgium-20}};
            \addplot[fill] table [x, y, col sep=comma] {data/hist-COLLECTION.n-k-15-FORMULATION.DWP-NETWORK.Belgium-20.csv};
            \node[fill=white,anchor=north east] at (rel axis cs:0.98,0.97) {\texttt{Belgium-20}};
            \addplot[fill=white] table [x, y, col sep=comma] {data/hist-COLLECTION.n-k-15-FORMULATION.CRDWP-NETWORK.Belgium-20.csv};
        \end{axis}
    \end{tikzpicture}
\end{subfigure}
\begin{subfigure}[b]{0.32\textwidth}
    \begin{tikzpicture}
        \begin{axis}[height=4.0cm,ybar,enlargelimits=false,bar width=2.00pt,
                     xmajorgrids={false},xmin=-5.0,xmax=115.0,ymin=0.0,
                     ymax=55.0,xticklabels={,,},tick pos=left,xtick={0,20,40,60,80,100},
                     yticklabels={,,}]
            \node[fill=white,anchor=north east] at (rel axis cs:0.98,0.97) {\texttt{GasLib-24}};
            \addplot[fill] table [x, y, col sep=comma] {data/hist-COLLECTION.n-k-15-FORMULATION.DWP-NETWORK.GasLib-24.csv};
            \node[fill=white,anchor=north east] at (rel axis cs:0.98,0.97) {\texttt{GasLib-24}};
            \addplot[fill=white] table [x, y, col sep=comma] {data/hist-COLLECTION.n-k-15-FORMULATION.CRDWP-NETWORK.GasLib-24.csv};
        \end{axis}
    \end{tikzpicture}
\end{subfigure} \\
\begin{subfigure}[b]{0.345\textwidth}
    \begin{tikzpicture}
        \begin{axis}[height=4.0cm,ybar,enlargelimits=false,bar width=2.00pt,
                     ylabel=Scenarios Solved ($\%$),xticklabels={,,},ymin=0.0,
                     xmajorgrids={false},xmin=-5.0,xmax=115.0,xtick={0,20,40,60,80,100},
                     tick pos=left,ymax=55.0]
            \node[fill=white,anchor=north east] at (rel axis cs:0.98,0.97) {\texttt{GasLib-40}};
            \addplot[fill] table [x, y, col sep=comma] {data/hist-COLLECTION.n-k-15-FORMULATION.DWP-NETWORK.GasLib-40.csv};
            \node[fill=white,anchor=north east] at (rel axis cs:0.98,0.97) {\texttt{GasLib-40}};
            \addplot[fill=white] table [x, y, col sep=comma] {data/hist-COLLECTION.n-k-15-FORMULATION.CRDWP-NETWORK.GasLib-40.csv};
        \end{axis}
    \end{tikzpicture}
\end{subfigure}
\begin{subfigure}[b]{0.32\textwidth}
    \hspace{0.5em}
    \begin{tikzpicture}
        \begin{axis}[height=4.0cm,ybar,enlargelimits=false,bar width=2.00pt,
                     xmajorgrids={false},xmin=-5.0,xmax=115.0,ymin=0.0,
                     ymax=55.0,xticklabels={,,},tick pos=left,xtick={0,20,40,60,80,100},
                     yticklabels={,,}]
            \node[fill=white,anchor=north east] at (rel axis cs:0.98,0.97) {\texttt{GasLib-134}};
            \addplot[fill] table [x, y, col sep=comma] {data/hist-COLLECTION.n-k-15-FORMULATION.DWP-NETWORK.GasLib-134.csv};
            \node[fill=white,anchor=north east] at (rel axis cs:0.98,0.97) {\texttt{GasLib-134}};
            \addplot[fill=white] table [x, y, col sep=comma] {data/hist-COLLECTION.n-k-15-FORMULATION.CRDWP-NETWORK.GasLib-134.csv};
        \end{axis}
    \end{tikzpicture}
\end{subfigure}
\begin{subfigure}[b]{0.32\textwidth}
    \begin{tikzpicture}
        \begin{axis}[height=4.0cm,ybar,enlargelimits=false,bar width=2.00pt,
                     xmajorgrids={false},xmin=-5.0,xmax=115.0,ymin=0.0,
                     ymax=55.0,xticklabels={,,},tick pos=left,xtick={0,20,40,60,80,100},
                     yticklabels={,,}]
            \node[fill=white,anchor=north east] at (rel axis cs:0.98,0.97) {\texttt{GasLib-135}};
            \addplot[fill] table [x, y, col sep=comma] {data/hist-COLLECTION.n-k-15-FORMULATION.DWP-NETWORK.GasLib-135.csv};
            \node[fill=white,anchor=north east] at (rel axis cs:0.98,0.97) {\texttt{GasLib-135}};
            \addplot[fill=white] table [x, y, col sep=comma] {data/hist-COLLECTION.n-k-15-FORMULATION.CRDWP-NETWORK.GasLib-135.csv};
        \end{axis}
    \end{tikzpicture}
\end{subfigure} \\
\begin{subfigure}[b]{0.345\textwidth}
    \begin{tikzpicture}
        \begin{axis}[height=4.0cm,ybar,enlargelimits=false,bar width=2.00pt,
                     ylabel=Scenarios Solved ($\%$),xtick={0,20,40,60,80,100},
                     xmajorgrids={false},xmin=-5.0,xmax=115.0,ymin=0.0,
                     xlabel={Load Delivered ($\%$)},
                     tick pos=left,ymax=55.0]
            \node[anchor=north east] at (rel axis cs:0.98,0.97) {\texttt{NA-154}};
            \addplot[fill] table [x, y, col sep=comma] {data/hist-COLLECTION.n-k-15-FORMULATION.DWP-NETWORK.NA-154.csv};
            \node[anchor=north east] at (rel axis cs:0.98,0.97) {\texttt{NA-154}};
            \addplot[fill=white] table [x, y, col sep=comma] {data/hist-COLLECTION.n-k-15-FORMULATION.CRDWP-NETWORK.NA-154.csv};
        \end{axis}
    \end{tikzpicture}
\end{subfigure}
\begin{subfigure}[b]{0.32\textwidth}
    \hspace{0.5em}
    \begin{tikzpicture}
        \begin{axis}[height=4.0cm,ybar,enlargelimits=false,bar width=2.00pt,
                     xmajorgrids={false},xmin=-5.0,xmax=115.0,ymin=0.0,
                     ymax=55.0,tick pos=left,xtick={0,20,40,60,80,100},
                     xlabel={Load Delivered ($\%$)},
                     yticklabels={,,}]
            \node[fill=white,anchor=north east] at (rel axis cs:0.98,0.97) {\texttt{GasLib-582}};
            \addplot[fill] table [x, y, col sep=comma] {data/hist-COLLECTION.n-k-15-FORMULATION.DWP-NETWORK.GasLib-582.csv};
            \node[fill=white,anchor=north east] at (rel axis cs:0.98,0.97) {\texttt{GasLib-582}};
            \addplot[fill=white] table [x, y, col sep=comma] {data/hist-COLLECTION.n-k-15-FORMULATION.CRDWP-NETWORK.GasLib-582.csv};
        \end{axis}
    \end{tikzpicture}
\end{subfigure}
\begin{subfigure}[b]{0.32\textwidth}
    \begin{tikzpicture}
        \begin{axis}[height=4.0cm,ybar,enlargelimits=false,bar width=2.00pt,
                     xmajorgrids={false},xmin=-5.0,xmax=115.0,ymin=0.0,
                     ymax=55.0,tick pos=left,xtick={0,20,40,60,80,100},
                     xlabel={Load Delivered ($\%$)},
                     yticklabels={,,}]
            \node[fill=white,anchor=north east] at (rel axis cs:0.98,0.97) {\texttt{GasLib-4197}};
            \addplot[fill] table [x, y, col sep=comma] {data/hist-COLLECTION.n-k-15-FORMULATION.DWP-NETWORK.GasLib-4197.csv};
            \node[fill=white,anchor=north east] at (rel axis cs:0.98,0.97) {\texttt{GasLib-4197}};
            \addplot[fill=white] table [x, y, col sep=comma] {data/hist-COLLECTION.n-k-15-FORMULATION.CRDWP-NETWORK.GasLib-4197.csv};
        \end{axis}
    \end{tikzpicture}
\end{subfigure}
    \caption{Histograms of load delivered over randomized $N{-}k$ multi-contingency scenarios using the \eqref{eqn:minqp} and \eqref{eqn:micqp} formulations, where $k$ corresponds to $15\%$ of all edges. Note that each pair of histograms summarizes instances solved by \emph{both} \eqref{eqn:minqp} and \eqref{eqn:micqp}.}
    \label{figure:n-k-15-histogram}
\end{figure}
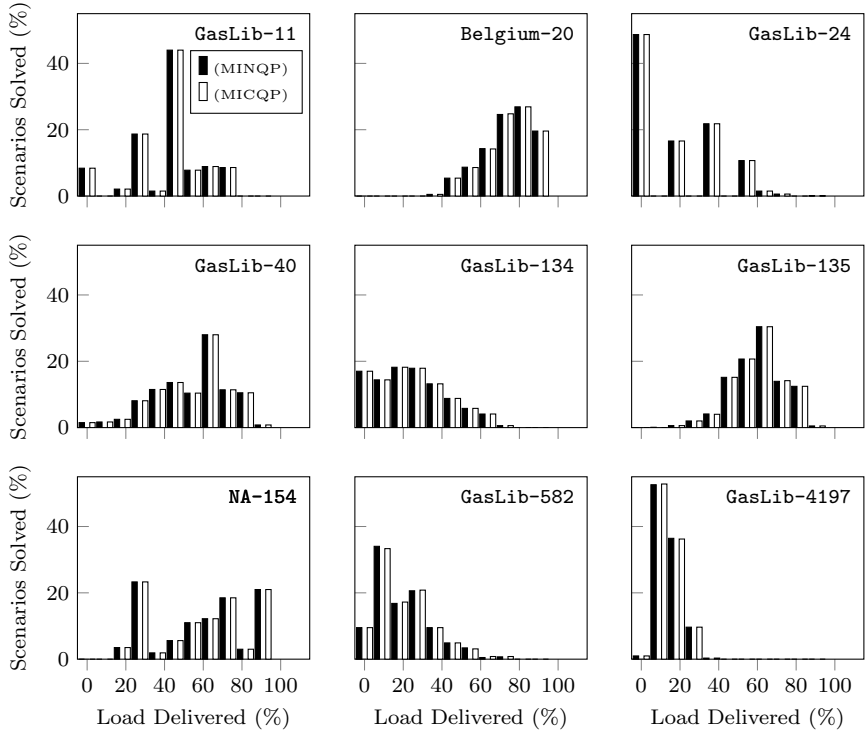

Figure \ref{figure:n-k-15-histogram} displays nine histograms that compare the proportion of load delivered across \emph{solved} damage scenarios for each network while using the two formulations.
Most importantly, these histograms display the similarity of the results achieved while using the relaxation-based formulation.
These results also indicate substantial qualitative differences in the hypothetical robustness of each network.
For example, larger networks like \texttt{GasLib-582} and \texttt{GasLib-4197} are shown to be highly sensitive to the $15\%$ damage scenarios, where often only $10\%$ to $30\%$ of load can be delivered.
The \texttt{Belgium-20} network appears less vulnerable and is often capable of serving between $70\%$ and $100\%$ of the load under severe contingencies.
Finally, for some smaller networks (e.g., \texttt{GasLib-11}, \texttt{GasLib-24}, \texttt{GasLib-134}), many scenarios result in zero or nearly zero deliverable load.
This simple analysis shows the utility of the MLD problem for understanding broad characteristics of gas network robustness.

\subsection{Computational Performance}
\label{subsection:computational-performance}
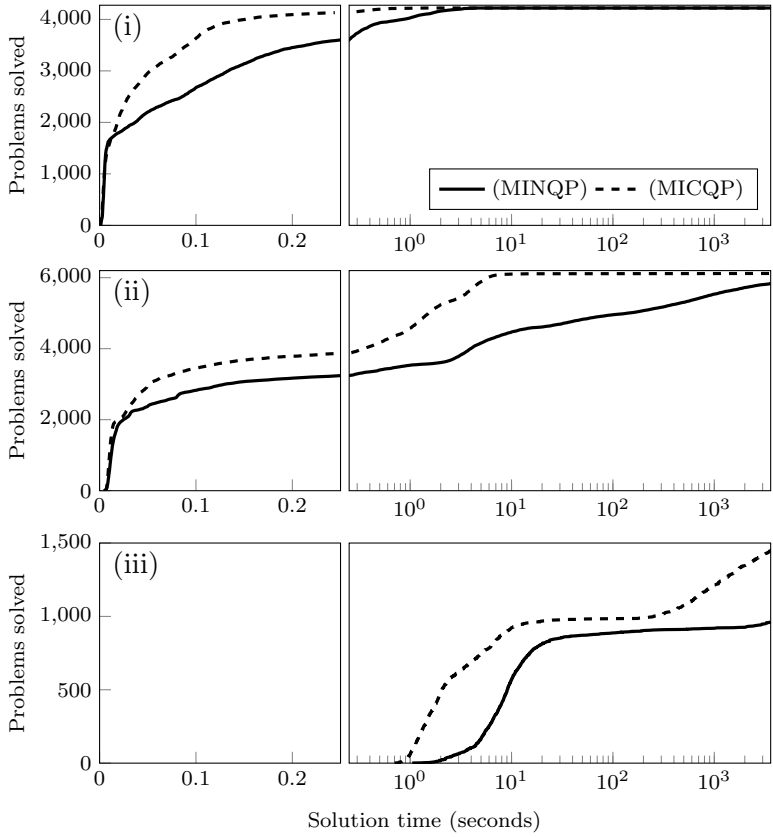
\begin{figure}[t]

\centering
\begin{subfigure}[b]{1.0\linewidth}
    \centering
    \begin{tikzpicture}[baseline]
        \begin{axis}[legend cell align=left,enlargelimits=false,
                     tick pos=left,height=4.5cm,ylabel=Problems solved,
                     scaled x ticks=false,xticklabel style={/pgf/number format/fixed},
                     width=0.4\linewidth,ymin=0,ymax=4275,restrict x to domain=0:0.250,xmin=0.0,xmax=0.250]
            \pgfplotstableread[col sep = comma]{data/performance-tens-FORMULATION.DWP.csv}{\dwp};
            \addplot[very thick] table [x = solve_time, y = num_problems_solved]{\dwp};
            \pgfplotstableread[col sep = comma]{data/performance-tens-FORMULATION.CRDWP.csv}{\crdwp};
            \addplot[very thick, dashed] table [x = solve_time, y = num_problems_solved]{\crdwp};
            \node[anchor=north west] at (rel axis cs:0.02,1.0) {\large{(i)}};
        \end{axis}
    \end{tikzpicture}
    \begin{tikzpicture}[baseline]
        \begin{semilogxaxis}[legend cell align=left,enlargelimits=false,
                             tick pos=left,height=4.5cm,
                             ytick=\empty,width=0.6\linewidth,scaled ticks=false,
                             xlabel absolute,legend style={at={(0.97, 0.07)},anchor=south east},
                             legend columns=2,ymin=0,ymax=4275,xmin=0.25,xmax=3600.0]
            \pgfplotstableread[col sep = comma]{data/performance-tens-FORMULATION.DWP.csv}{\dwp};
            \addplot[very thick] table [x = solve_time, y = num_problems_solved]{\dwp};
            \addlegendentry{\eqref{eqn:minqp}}
            \pgfplotstableread[col sep = comma]{data/performance-tens-FORMULATION.CRDWP.csv}{\crdwp};
            \addplot[very thick, dashed] table [x = solve_time, y = num_problems_solved]{\crdwp};
            \addlegendentry{\eqref{eqn:micqp}}
        \end{semilogxaxis}
    \end{tikzpicture}
\end{subfigure} \\
\begin{subfigure}[b]{1.0\linewidth}
    \centering
    \begin{tikzpicture}[baseline]
        \begin{axis}[legend cell align=left,enlargelimits=false,
                     tick pos=left,height=4.5cm,ylabel=Problems solved,
                     scaled x ticks=false,xticklabel style={/pgf/number format/fixed},
                     width=0.4\linewidth,ymin=0,ymax=6200,restrict x to domain=0:0.25,xmin=0.0,xmax=0.25]
            \pgfplotstableread[col sep = comma]{data/performance-hundreds-FORMULATION.DWP.csv}{\dwp};
            \addplot[very thick] table [x = solve_time, y = num_problems_solved]{\dwp};
            \pgfplotstableread[col sep = comma]{data/performance-hundreds-FORMULATION.CRDWP.csv}{\crdwp};
            \addplot[very thick, dashed] table [x = solve_time, y = num_problems_solved]{\crdwp};
            \node[anchor=north west] at (rel axis cs:0.02,1.00) {\large{(ii)}};
        \end{axis}
    \end{tikzpicture}
    \begin{tikzpicture}[baseline]
        \begin{semilogxaxis}[legend cell align=left,enlargelimits=false,
                             tick pos=left,height=4.5cm,
                             ytick=\empty,width=0.6\linewidth,
                             xlabel absolute,legend style={at={(0.97, 0.04)},anchor=south east},
                             legend columns=2,ymin=0,ymax=6200,xmin=0.25,xmax=3600.0]
            \pgfplotstableread[col sep = comma]{data/performance-hundreds-FORMULATION.DWP.csv}{\dwp};
            \addplot[very thick] table [x = solve_time, y = num_problems_solved]{\dwp};
            \pgfplotstableread[col sep = comma]{data/performance-hundreds-FORMULATION.CRDWP.csv}{\crdwp};
            \addplot[very thick, dashed] table [x = solve_time, y = num_problems_solved]{\crdwp};
        \end{semilogxaxis}
    \end{tikzpicture}
\end{subfigure} \\
\begin{subfigure}[b]{1.0\linewidth}
    \centering
    \begin{tikzpicture}[baseline]
        \begin{axis}[legend cell align=left,enlargelimits=false,
                     tick pos=left,height=4.5cm,ylabel=Problems solved,
                     scaled x ticks=false,xticklabel style={/pgf/number format/fixed},
                     width=0.4\linewidth,ymin=0,ymax=1500,restrict x to domain=0:0.25,xmin=0.0,xmax=0.25]
            \pgfplotstableread[col sep = comma]{data/performance-thousands-FORMULATION.DWP.csv}{\dwp};
            \addplot[very thick] table [x = solve_time, y = num_problems_solved]{\dwp};
            \pgfplotstableread[col sep = comma]{data/performance-thousands-FORMULATION.CRDWP.csv}{\crdwp};
            \addplot[very thick, dashed] table [x = solve_time, y = num_problems_solved]{\crdwp};
            \node[anchor=north west] at (rel axis cs:0.02,1.00) {\large{(iii)}};
        \end{axis}
    \end{tikzpicture}
    \begin{tikzpicture}[baseline]
        \begin{semilogxaxis}[legend cell align=left,enlargelimits=false,
                             tick pos=left,height=4.5cm,x tick label style={/pgf/number format/fixed},
                             xlabel={\hspace{-0.3\linewidth}Solution time (seconds)},
                             ytick=\empty,width=0.6\linewidth,
                             xlabel absolute,legend style={at={(0.97, 0.04)},anchor=south east},
                             legend columns=2,ymin=0,ymax=1500,xmin=0.25,xmax=3600.0]
            \pgfplotstableread[col sep = comma]{data/performance-thousands-FORMULATION.DWP.csv}{\dwp};
            \addplot[very thick] table [x = solve_time, y = num_problems_solved]{\dwp};
            \pgfplotstableread[col sep = comma]{data/performance-thousands-FORMULATION.CRDWP.csv}{\crdwp};
            \addplot[very thick, dashed] table [x = solve_time, y = num_problems_solved]{\crdwp};
        \end{semilogxaxis}
    \end{tikzpicture}
\end{subfigure}
    \caption{Performance profiles comparing the efficiency of \eqref{eqn:minqp} and \eqref{eqn:micqp} formulations over the instances described in Sections \ref{subsection:n-1_experiments} and \ref{subsection:n-k_experiments}. Here, the performance profiles are divided into three categories for (i) networks containing tens of nodes; (ii) networks containing hundreds of nodes; and (iii) networks containing thousands of nodes (i.e., \texttt{GasLib-4197}).}
    \label{figure:performance}
\end{figure}

This section compares the performance of \eqref{eqn:minqp} and \eqref{eqn:micqp} MLD formulations using the benchmark instances described in Sections \ref{subsection:n-1_experiments} and \ref{subsection:n-k_experiments}.
The performance profiles for these instances are shown in Figure \ref{figure:performance} and divided into three categories: (i) networks containing tens of nodes; (ii) networks containing hundreds of nodes; and (iii) networks containing thousands of nodes (i.e., \texttt{GasLib-4197}).
In all such categories, it is shown that the \eqref{eqn:micqp} formulation is able to solve substantially greater numbers of problems than \eqref{eqn:minqp} in shorter amounts of time.
For networks with tens of nodes, both formulations are capable of solving all instances in less than ten seconds.
For networks with hundreds of nodes, \eqref{eqn:micqp} is capable of solving most instances within ten seconds, while \eqref{eqn:minqp} requires hundreds or thousands of seconds.
For networks with thousands of nodes, \eqref{eqn:micqp} solves a much greater number of instances within the alotted hour time limit, although solve times are much longer than for problems containing networks of tens or hundreds of nodes.

\subsection{Synthetic Earthquake Damage Scenarios}
\label{subsection:earthquake_experiments}
This subsection provides a proof of concept application to demonstrate the use of the MLD problem in the context of risk assessment for deterministic and uncertain spatial hazards.
In each scenario, damage to a network is assumed to be caused by an earthquake with a fixed magnitude and epicenter.
For an earthquake, the probability of damage to a pipeline component is commonly represented as a function of peak ground acceleration (PGA) and peak ground velocity (PGV).
These relationships for PGA and PGV are typically expressed as functions of earthquake magnitude and distance from seismographic rupture.
In this paper, the relationships derived by \citet{campbell1997empirical} are used, although their details are omitted here.
All earthquakes are assumed to arise from strike-slip faulting at a depth of one kilometer in a region of firm soil.

Given the PGA and PGV at a component's point in space, the vulnerability of the component is then modeled using the fragility approach of \citet{lanzano2014seismic} for continuous pipelines in the presence of strong ground shaking.
Specifically, the probability of damage is computed as a function of PGV, and node-connecting components exceeding the damage threshold for the first risk state (``very limited loss'') are assumed to be nonfunctional.
Again, for the purpose of brevity, these relationships are omitted here.
For simplicity, we make three additional assumptions to the model: (i) only node-connecting components are affected by PGA and PGV; (ii) all node-connecting components are assumed structurally equivalent as ``continuous pipelines,'' and (iii) the distance from the epicenter to each node-connecting component is the minimum distance between the epicenter and the locations of the connecting junctions.

\begin{figure}[t]
    \input{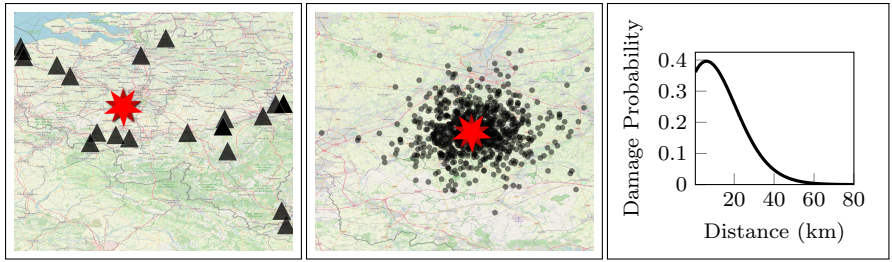}
    \caption{Illustrations of the earthquake scenarios. The first (left) shows the position of an earthquake epicenter (red star) for the \texttt{Belgium-20} network and its spatial relation to network junctions (black triangles). The second illustration (center) shows the placement of a mean earthquake epicenter (red star) for the \texttt{Belgium-20} network and random, normally distributed epicenters surrounding it (black circles). The last (right) shows a fragility curve derived from \citet{campbell1997empirical} and \citet{lanzano2014seismic} for a magnitude $8.0$ earthquake.}
    \label{figure:earthquake}
\end{figure}

The networks \texttt{Belgium-20}, \texttt{GasLib-40}, \texttt{GasLib-135}, \texttt{NA-154}, \texttt{GasLib-582}, and \texttt{GasLib-4197} are considered for earthquake damage scenarios, as the remaining three networks do not contain geolocation data for components.
Figure \ref{figure:earthquake} depicts three illustrations relevant to the parameterization of scenarios.

\paragraph{\textbf{Deterministic Earthquake, Stochastic Fragility Scenarios.}}
The first set of earthquake scenarios is intended to demonstrate the applicability of the MLD method when analyzing network vulnerability to a known (i.e., deterministic) natural hazard.
To this end, one earthquake is considered per network, each described by a fixed magnitude and epicenter.
This type of scenario is illustrated pictorially in the first image of Figure \ref{figure:earthquake}.
The local magnitude of each earthquake is assumed to be $8.0$, while each epicenter is assumed to be the center of a $k$-means cluster containing the largest number of junctions, where $k = 5$.
Then, using the PGV model of \citet{campbell1997empirical} and the probabilistic fragility approach of \citet{lanzano2014seismic}, one thousand damage scenarios are generated per network, where in each, the operational status per node-connecting component is determined via a uniform random sampling and comparison with the probability of damage.
An example fragility curve is depicted in the last image of Figure \ref{figure:earthquake}, which relates a component's distance from the epicenter of a magnitude $8.0$ earthquake to the probability of damage.


\begin{table}[t]
    \begin{center}
    \begin{tabular}{c|c|c|c|c|c|c|}
        \cline{2-7} & \multicolumn{3}{c|}{MINQP} & \multicolumn{3}{c|}{MICQP} \\
        \cline{1-7} \multicolumn{1}{|c|}{Network} & Opt. (\%) & Lim. (\%) & Inf. (\%) & Opt. (\%) & Lim. (\%) & Inf. (\%) \\
        \cline{1-7} \multicolumn{1}{|c|}{\texttt{Belgium-20}} & $100.0$ & $0.0$ & $0.0$ & $100.0$ & $0.0$ & $0.0$ \\
        \cline{1-7} \multicolumn{1}{|c|}{\texttt{GasLib-40}} & $100.0$ & $0.0$ & $0.0$ & $100.0$ & $0.0$ & $0.0$ \\
        \cline{1-7} \multicolumn{1}{|c|}{\texttt{GasLib-135}} & $98.3$ & $1.7$ & $0.0$ & $100.0$ & $0.0$ & $0.0$ \\
        \cline{1-7} \multicolumn{1}{|c|}{\texttt{NA-154}} & $100.0$ & $0.0$ & $0.0$ & $100.0$ & $0.0$ & $0.0$ \\
        \cline{1-7} \multicolumn{1}{|c|}{\texttt{GasLib-582}} & $58.7$ & $41.3$ & $0.0$ & $99.8$ & $0.2$ & $0.0$ \\
        \cline{1-7} \multicolumn{1}{|c|}{\texttt{GasLib-4197}} & $2.9$ & $95.3$ & $1.8$ & $86.0$ & $14.0$ & $0.0$ \\
        \cline{1-7}
    \end{tabular}
    \end{center}
    \caption{Comparison of solution statuses over the set of deterministic earthquake scenarios.}
    \label{table:earthquake-1}
\end{table}

Table \ref{table:earthquake-1} displays statistics of solver statuses across all damage scenarios for each network and formulation.
For all except the \texttt{GasLib-135}, \texttt{GasLib-582}, and \texttt{GasLib-4197} networks, optimal solutions to instances are found within the one hour time limit.
For \texttt{GasLib-135}, the \eqref{eqn:minqp} formulation is unable to prove optimality on $17$ instances.
For \texttt{GasLib-582} and \texttt{GasLib-4197}, a larger proportion of \eqref{eqn:minqp} instances cannot be solved within the time limit.
Even for the convex relaxation, two instances cannot be solved for the \texttt{GasLib-582} network.
Nonetheless, comparing the \eqref{eqn:minqp} and \eqref{eqn:micqp} formulations shows the benefit of the relaxation-based approach, which is capable of solving a larger proportion of instances.
We note that $96$ of the thousand \texttt{GasLib-4197} \eqref{eqn:minqp} instances are claimed to be infeasible using default \textsc{Gurobi} parameters.
This number is reduced to $18$ using \texttt{NumericFocus=3}.

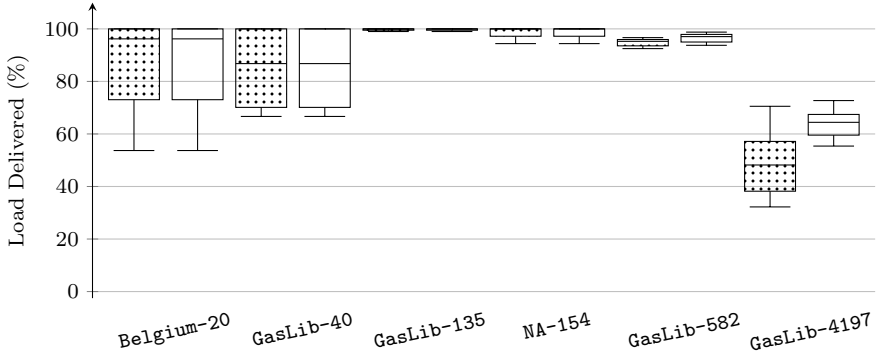
\begin{figure}[t]

\begin{tikzpicture}
    \pgfplotstableread[col sep=comma]{data/earthquake_1_boxplots.csv}\csvdata
    \pgfplotstabletranspose\datatransposed{\csvdata} 
    \begin{axis}[boxplot/draw direction = y, x axis line style = {opacity=0},
                 enlarge x limits={abs=6pt},
                 axis x line* = bottom, axis y line = left, height = 5.5cm,
                 ymajorgrids, width=\textwidth, xtick = {1, 2, 3, 4, 5, 6, 7, 8, 9, 10, 11, 12},
                 xticklabel style = {align=center, font=\small, rotate=12, xshift=-0.30cm},
                 xticklabels = {,\texttt{Belgium-20},,
                 \texttt{GasLib-40},,\texttt{GasLib-135},,\texttt{NA-154},,
                 \texttt{GasLib-582},,\texttt{GasLib-4197}}, xtick style = {draw=none},
                 ylabel = {Load Delivered ($\%$)}, ymin = -2.5, ymax = 110.0,
                 ytick = {0.0, 20.0, 40.0, 60.0, 80.0, 100.0}]
        \addplot+[boxplot, pattern = dots, draw=black] table[y index=1] {\datatransposed};
        \addplot+[boxplot, draw=black] table[y index=2] {\datatransposed};
        \addplot+[boxplot, pattern = dots, draw=black] table[y index=3] {\datatransposed};
        \addplot+[boxplot, draw=black] table[y index=4] {\datatransposed};
        \addplot+[boxplot, pattern = dots, draw=black] table[y index=5] {\datatransposed};
        \addplot+[boxplot, draw=black] table[y index=6] {\datatransposed};
        \addplot+[boxplot, pattern = dots, draw=black] table[y index=7] {\datatransposed};
        \addplot+[boxplot, draw=black] table[y index=8] {\datatransposed};
        \addplot+[boxplot, pattern = dots, draw=black] table[y index=9] {\datatransposed};
        \addplot+[boxplot, draw=black] table[y index=10] {\datatransposed};
        \addplot+[boxplot, pattern = dots, draw=black] table[y index=11] {\datatransposed};
        \addplot+[boxplot, draw=black] table[y index=12] {\datatransposed};
    \end{axis}
\end{tikzpicture}
    \caption{Boxplots of load delivered over \emph{solved} earthquake damage scenarios per network, where the magnitude and epicenter for each network are assumed to be fixed. Here, dotted and clear boxplots correspond to \eqref{eqn:minqp} and \eqref{eqn:micqp} formulations, respectively. Note that each pair of boxplots summarizes only instances solved by \emph{both} \eqref{eqn:minqp} and \eqref{eqn:micqp}.}
    \label{figure:earthquake_1_boxplots}
\end{figure}

Figure \ref{figure:earthquake_1_boxplots} displays boxplots comparing the maximal proportion of load delivered across \emph{solved} damage scenarios for each network and formulation.
Here, results obtained using the \eqref{eqn:minqp} and \eqref{eqn:micqp} formulations of the MLD problem are shown to be remarkably similar.
This demonstrates the utility of the relaxation-based approach, which provides outcomes comparable to the mixed-integer nonconvex formulation at a smaller computational cost.
The boxplots also show substantial qualitative differences in the (hypothetical) vulnerability among networks.
For example, the \texttt{GasLib-40} network is shown to have great variability in maximal load delivered, while \texttt{GasLib-135}, \texttt{NA-154}, and \texttt{GasLib-582} are mostly unaffected by the hypothetical hazard.
Additionally, some networks (e.g., \texttt{NA-154}, \texttt{GasLib-582}) predict small ranges of load delivered, while others (e.g., \texttt{Belgium-20}) appear to carry greater uncertainty.

The large discrepancies in boxplots for the \texttt{GasLib-4197} damage scenarios likely originate from two sources.
First, the pair of boxplots is representative of only the subset of instances solved by \emph{both} formulations, which is relatively small.
Second, however, are the mathematical differences in nonconvex and relaxed formulations.
Notably, the minima across the solvable instances differ by around twenty percent.
This could be a consequence of the relaxation, which theoretically predicts maximal load values greater than or equal to the nonconvex formulation.
Compared to other networks, these larger discrepancies in predicted deliverable load could be a manifestation of the component relaxations, whose errors are further aggregated as the network size grows.

\paragraph{\textbf{Stochastic Earthquake, Stochastic Fragility Scenarios.}}
The second set of earthquake scenarios is intended to demonstrate the applicability of the MLD method when analyzing network vulnerability to a stochastic (i.e., uncertain) natural hazard.
To this end, multiple earthquakes are considered per network, where each is randomly sampled assuming normally distributed magnitudes and epicenters.
This situation is illustrated in the second image of Figure \ref{figure:earthquake}.
Here, the mean local magnitude of each earthquake is assumed to be $8.0$ with a standard deviation of $0.25$, and each mean epicenter is again assumed to be the center of a $k$-means cluster containing the largest number of junctions, where $k = 5$, and where a distance-based standard deviation of ten kilometers is assumed.
Using the models of \citet{campbell1997empirical} and \citet{lanzano2014seismic}, one thousand random earthquake scenarios are generated per network, where in each, the status per node-connecting component is again determined via a uniform random sampling and comparison with the probability of damage.


\begin{table}[t]
    \begin{center}
    \begin{tabular}{c|c|c|c|c|c|c|}
        \cline{2-7} & \multicolumn{3}{c|}{MINQP} & \multicolumn{3}{c|}{MICQP} \\
        \cline{1-7} \multicolumn{1}{|c|}{Network} & Opt. (\%) & Lim. (\%) & Inf. (\%) & Opt. (\%) & Lim. (\%) & Inf. (\%) \\
        \cline{1-7} \multicolumn{1}{|c|}{\texttt{Belgium-20}} & $100.0$ & $0.0$ & $0.0$ & $100.0$ & $0.0$ & $0.0$ \\
        \cline{1-7} \multicolumn{1}{|c|}{\texttt{GasLib-40}} & $100.0$ & $0.0$ & $0.0$ & $100.0$ & $0.0$ & $0.0$ \\
        \cline{1-7} \multicolumn{1}{|c|}{\texttt{GasLib-135}} & $97.3$ & $2.7$ & $0.0$ & $100.0$ & $0.0$ & $0.0$ \\
        \cline{1-7} \multicolumn{1}{|c|}{\texttt{NA-154}} & $100.0$ & $0.0$ & $0.0$ & $100.0$ & $0.0$ & $0.0$ \\
        \cline{1-7} \multicolumn{1}{|c|}{\texttt{GasLib-582}} & $97.1$ & $2.9$ & $0.0$ & $100.0$ & $0.0$ & $0.0$ \\
        \cline{1-7} \multicolumn{1}{|c|}{\texttt{GasLib-4197}} & $4.4$ & $89.7$ & $5.9$ & $82.8$ & $17.2$ & $0.0$ \\
        \cline{1-7}
    \end{tabular}
    \end{center}
    \caption{Comparison of termination statuses over the second set of earthquake scenarios.}
    \label{table:earthquake-2}
\end{table}

Table \ref{table:earthquake-2} repeats the format of Table \ref{table:earthquake-1} to display aggregate statistics of solver termination statuses across all stochastic earthquake scenarios.
For most scenarios, globally optimal solutions to MLD instances are found within the prescribed one hour time limit.
For scenarios based on the \texttt{GasLib-135}, \texttt{GasLib-582}, and \texttt{GasLib-4197} networks, however, some instances cannot be solved.
As in the deterministic scenario analysis, comparing the \eqref{eqn:minqp} and \eqref{eqn:micqp} formulations shows the benefit of the relaxation-based approach.
Furthermore, $133$ and one \texttt{GasLib-4197} \eqref{eqn:minqp} and \eqref{eqn:micqp} instances, respectively, are claimed to be infeasible using default \textsc{Gurobi} parameters.
These numbers are reduced to $59$ and zero, respectively, using \texttt{NumericFocus=3}.

\begin{figure}[t]

\begin{tikzpicture}
    \pgfplotstableread[col sep=comma]{data/earthquake_2_boxplots.csv}\csvdata
    \pgfplotstabletranspose\datatransposed{\csvdata} 
    \begin{axis}[boxplot/draw direction = y, x axis line style = {opacity=0},
                 enlarge x limits={abs=6pt},
                 axis x line* = bottom, axis y line = left, height = 5.5cm,
                 ymajorgrids, width=\textwidth, xtick = {1, 2, 3, 4, 5, 6, 7, 8, 9, 10, 11, 12},
                 xticklabel style = {align=center, font=\small, rotate=12, xshift=-0.30cm},
                 xticklabels = {,\texttt{Belgium-20},,
                 \texttt{GasLib-40},,\texttt{GasLib-135},,\texttt{NA-154},,
                 \texttt{GasLib-582},,\texttt{GasLib-4197}}, xtick style = {draw=none},
                 ylabel = {Load Delivered ($\%$)}, ymin = -2.5, ymax = 110.0,
                 ytick = {0.0, 20.0, 40.0, 60.0, 80.0, 100.0}]
        \addplot+[boxplot, pattern = dots, draw=black] table[y index=1] {\datatransposed};
        \addplot+[boxplot, draw=black] table[y index=2] {\datatransposed};
        \addplot+[boxplot, pattern = dots, draw=black] table[y index=3] {\datatransposed};
        \addplot+[boxplot, draw=black] table[y index=4] {\datatransposed};
        \addplot+[boxplot, pattern = dots, draw=black] table[y index=5] {\datatransposed};
        \addplot+[boxplot, draw=black] table[y index=6] {\datatransposed};
        \addplot+[boxplot, pattern = dots, draw=black] table[y index=7] {\datatransposed};
        \addplot+[boxplot, draw=black] table[y index=8] {\datatransposed};
        \addplot+[boxplot, pattern = dots, draw=black] table[y index=9] {\datatransposed};
        \addplot+[boxplot, draw=black] table[y index=10] {\datatransposed};
        \addplot+[boxplot, pattern = dots, draw=black] table[y index=11] {\datatransposed};
        \addplot+[boxplot, draw=black] table[y index=12] {\datatransposed};
    \end{axis}
\end{tikzpicture}
    \caption{Boxplots of load delivered over earthquake damage scenarios per network, where both the magnitude and epicenter for each scenario are normally distributed. Here, dotted and clear boxplots correspond to \eqref{eqn:minqp} and \eqref{eqn:micqp} formulations, respectively. Note that each pair of boxplots summarizes only instances solved by \emph{both} \eqref{eqn:minqp} and \eqref{eqn:micqp}.}
    \label{figure:earthquake_2_boxplots}
\end{figure}
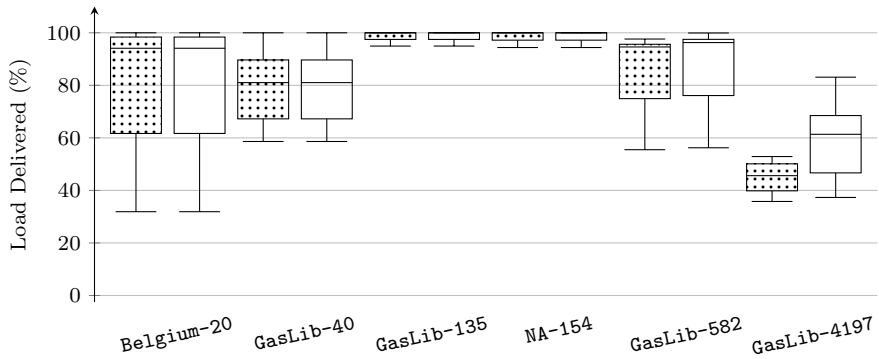

Figure \ref{figure:earthquake_2_boxplots} displays boxplots comparing the proportion of load delivered across all scenarios.
Again, the \eqref{eqn:minqp} and \eqref{eqn:micqp} formulations provide results that are qualitatively similar.
The boxplots also show similar differences in vulnerability as those observed in Figure \ref{figure:earthquake_2_boxplots}.
However, the greater variation in epicenters and magnitudes results in greater variation of the effects.

\section{Concluding Remarks}
\label{section:concluding_remarks}
This study introduces the MLD problem for natural gas transmission networks, which seeks to determine a feasible operating point for a damaged gas network while ensuring the maximal delivery of prioritized load.
This task is presented as three successive mathematical programming formulations.
First, a mixed-integer nonconvex program is formulated that embeds all physical and engineering requirements for operational feasibility.
Second, the first nonconvex program is reformulated \emph{exactly} as a mixed-integer nonconvex quadratic program by introducing new variables.
Finally, convex relaxations are applied to all nonconvex relationships in the former problem that involve pressures and mass flows, which results in a mixed-integer convex relaxed formulation.

To compare the efficacy of the second and third formulations, a rigorous benchmarking study is conducted over a large number of randomized multi-contingency scenarios on nine networks ranging in size from $11$ to $4197$ junctions.
First, MLD experiments based on $N{-}1$ (or single contingency) damage scenarios are conducted.
Second, $N{-}k$ (multi-contingency) damage experiments are performed in order to understand MLD tractability for multimodal failures.
A performance comparison of the \eqref{eqn:minqp} and \eqref{eqn:micqp} formulations is then conducted using results from the $N{-}1$ and $N{-}k$ experiments.
Finally, a proof of concept application based on network damage from a set of synthetically generated earthquakes demonstrates the application of the MLD problem to risk assessment for deterministic and stochastic spatial hazards.

These experimental results lead to three key conclusions.
First, the relaxed formulation \eqref{eqn:micqp} provides good bounds on the maximal deliverable load obtained from the full mixed-integer nonconvex formulation, \eqref{eqn:minqp}.
Second, the relaxation-based formulation is more computationally robust than the mixed-integer nonconvex formulation, and can solve larger proportions of challenging instances in much shorter amounts of time.
Finally, for the largest network (i.e., \texttt{GasLib-4197}), the relaxation-based approach begins to show its limitations.
For some challenging scenario types (e.g., $N{-}1$), large numbers of instances cannot be solved because of the relatively larger size of the network.

There are several potential studies that could extend the approaches developed in this study.
First, additional relaxation-based methods could be developed to more accurately and efficiently solve the MLD problem for gas networks containing thousands of nodes.
To facilitate this, another useful contribution could be the development of new benchmark instances whose sizes range between the sizes of \texttt{GasLib-582} and \texttt{GasLib-4197}.
Finally, the origin of numerical instabilities associated with the challenging \texttt{GasLib-4197} network, and especially the sources of claimed infeasibility by \textsc{Gurobi} for some \eqref{eqn:minqp} instances, should be thoroughly investigated.
This could involve developing new methods for processing network data or normalizing constraints.

\section*{Acknowledgments}
The authors express gratitude to Mary Ewers and David Fobes of Los Alamos National Laboratory for their construction of the \texttt{NA-154} network model.

\section*{Declarations}
\paragraph{\textbf{Funding}}
The Department of Homeland Security sponsored the production of this material under the U.S. Department of Energy contract for the management and operation of Los Alamos National Laboratory.
Los Alamos National Laboratory is operated by Triad National Security, LLC, for the National Nuclear Security Administration of the U.S. Department of Energy (Contract No. 89233218CNA000001).
By acceptance of this article, the publisher recognizes that the U.S. Government retains a nonexclusive, royalty-free license to publish or reproduce the published form of this contribution, or to allow others to do so, for U.S. Government purposes.
Los Alamos National Laboratory requests that the publisher identify this article as work performed under the auspices of the U.S. Department of Energy.

\paragraph{\textbf{Conflicts of Interest/Competing Interests}}
The authors have no conflicts of interest nor any competing interests.

\paragraph{\textbf{Availability of Data and Material}}
The \texttt{Belgium-20} and \texttt{NA-154} network data are available upon request, subject to approval.
The \texttt{GasLib} network data are publicly available at \href{https://gaslib.zib.de}{https://gaslib.zib.de}.

\paragraph{\textbf{Code Availability}}
The formulations in this paper are implemented in \textsc{GasModels}, publicly available at \href{https://github.com/lanl-ansi/GasModels.jl}{https://github.com/lanl-ansi/GasModels.jl}.
The source code for damage scenario preparation and analysis is available upon request, subject to approval.

\bibliographystyle{spbasic}
\bibliography{bibliography}
\end{document}